\documentclass[aos,MSNbibl,nameyear,seceqn,dvips]{arximspdf}
\usepackage{mathbh}

\doi{10.1214/14-AOS1274} 
\volume{43}
\issue{1}
\pubyear{2015}
\firstpage{215}
\lastpage{237}
\docsubty{FLA}

\makeatletter

\newproclaim{remark}{Remark}

\newcommand{\rrvert}{\vert}

\newcommand{\rrVert}{\Vert}
\newcommand{\llvert}{\vert}
\newcommand{\llVert}{\Vert}
\newcommand{\eqref}[1]{(\ref{#1})}

\newtheorem{theorem}{Theorem}[section]
\newtheorem{lemma}[theorem]{Lemma}
\newproclaim{example}[theorem]{Example}
\newtheorem{corollary}[theorem]{Corollary}
\newcommand{\pr}{\mathbb{P}}
\newcommand{\E}{\mathbb{E}}
\newcommand{\reals}{\mathbb{R}}
\newcommand{\ind}{\mathbh{1}}
\let\hat\widehat
\makeatother

\begin{document}
\begin{frontmatter}

\title{Consistency of spectral clustering in stochastic block models}
\runtitle{Consistency of spectral clustering}

\begin{aug}
\author[A]{\fnms{Jing}~\snm{Lei}\corref{}\thanksref{T1}\ead[label=e1]{jinglei@andrew.cmu.edu}}
and
\author[A]{\fnms{Alessandro}~\snm{Rinaldo}\thanksref{T2}\ead[label=e2]{arinaldo@cmu.edu}}
\runauthor{J. Lei and A. Rinaldo}
\thankstext{T1}{Supported by NSF Grant BCS-0941518, NSF Grant
DMS-14-07771 and NIH Grant MH057881.}
\thankstext{T2}{Supported by AFOSR and DARPA Grant FA9550-12-1-0392 and
NSF CAREER Grant DMS-11-49677.}

\affiliation{Carnegie Mellon University}

\address[A]{Department of Statistics\\
Carnegie Mellon University\\
132 Baker Hall/5000 Forbes Ave.\\
Pittsburgh, Pennsylvania 15213\\
USA\\
\printead{e1}\\
\phantom{E-mail:\ }\printead*{e2}}
\end{aug}

\received{\smonth{9} \syear{2014}}

%
\begin{abstract}
We analyze the performance of spectral clustering for community
extraction in stochastic
block models.
We show that, under mild conditions, spectral clustering
applied to the adjacency matrix of the network
can consistently recover hidden communities even when the order of the maximum
expected degree is as small as $\log n$, with $n$ the number of nodes.
This result applies to some popular
polynomial time spectral clustering algorithms and is further extended
to degree corrected
stochastic block
models using a spherical $k$-median spectral clustering method.
A key component of our analysis is a combinatorial bound on the
spectrum of binary random matrices, which is sharper than
the conventional matrix Bernstein inequality and may be of independent interest.
\end{abstract}
%

%
\end{frontmatter}

\section{Introduction}\label{sec1}
Network analysis is concerned with describing and modeling the joint occurrence
of random interactions among actors in a given population of interest.
In its simplest form, a network dataset over $n$ actors is a
simple undirected random graph on $n$ nodes, where the edges encode the realized
binary interactions among the nodes.
Examples include social networks (friendship between Facebook users,
blog following, twitter following, etc.),
biological networks (gene network, gene-protein network), information
network (email network, World Wide Web)
and many others.
A review of modeling and inference on network data
can be found in \citet{Kolaczyk09}, \citet{Newman09} and \citet
{Goldenberg10}.


Among the many existing statistical models for network data, the
stochastic block
model, henceforth SBM, of \citet{Holland83} stands out for both its simplicity
and expressive power.
In a SBM, the nodes are partitioned into $K < n$ disjoint groups, or
{\it
communities},
according to some latent random mechanism. Conditionally on the realized
but unobservable community
assignments, the edges then occur independently with probabilities
depending only on the
community membership of the nodes, so that nodes from the same
community will have
higher average degree of connectivity among themselves than compared to
the remaining
nodes (see Section~\ref{sec:setup} for details). Because of its simple
analytic form and its ability to capture the
emergence of communities, a feature commonly observed in
real network data, the SBM is certainly among the
most popular models for network data.

Within the SBM framework, the most important inferential task is
that of recovering the community membership of the nodes
from a \textit{single} observation of the network. To solve this problem,
in recent years researchers
have proposed a variety
procedures, which vary greatly
in their degrees of statistical accuracy and
computational complexity.
See, in particular, modularity maximization [\citet{NewmanG04}],
likelihood methods [\citet{BickelC09,ChoiWA12,ZhaoLZ12,AminiCBL12,Celisse12}],
method of moments [\citet{Anandkumar13}],
belief propagation [\citet{Decelle11}], convex optimization [\citet{ChenSX12}],
spectral clustering [\citet
{RoheCY11,siva.spectral:11,Jin12,Fishkind13,sarkar.bickel:13}] and its variants
[\citet{Coja-Oghlan10,ChaudhuriCT12}] and spectral embeddings
[\citeauthor{Sussman:12} (\citeyear
{Sussman:12}), \citet{Lyzinski13}].

Spectral clustering [see, e.g., \citet{VonLuxburg07}] is arguably one
of the most widely used
methods for community recovery. Broadly speaking, this procedure first performs
an eigen-decomposition of the adjacency matrix or the graph Laplacian.
Then the community membership is
inferred
by applying a clustering algorithm, typically $k$-means, to the (possibly
normalized) rows of the
matrix formed by the first few leading eigenvectors. Spectral
clustering is
easier to implement and computationally less demanding than many other
methods, most of which amount to computationally intractable combinatorial
searches. 
From a theoretical standpoint, spectral clustering has been shown to enjoy
good theoretical properties in denser stochastic block models where
the average degree grows faster
than $\log n$; see, for example, \citet{RoheCY11,Jin12,sarkar.bickel:13}.
In addition, spectral clustering has been empirically observed to
yield good performance even
in sparser regimes.
For example, it is recommended as the initial solution for a search
based procedure in
\citet{AminiCBL12}.
In computer science literature, spectral clustering is also a
standard procedure for graph partitioning and for
solving
the planted partition model, a special case of the SBM
[see, e.g., \citet{NgJW02}].

Despite its popularity and simplicity, the theoretical properties of
spectral clustering
are still not well understood in sparser SBM settings where the
magnitude of the maximum
expected node degree can be as small as $\log n$. This regime of
sparsity is in fact not covered by existing analyses of
the performance of spectral clustering for community recovery, which
postulate a denser network.
Indeed, \citet{RoheCY11,Fishkind13} require the expected node degree to
be almost linear
in $n$, while \citet{Jin12} requires polynomial growth. Analogous conditions
can be found elsewhere; see, for example,
\citeauthor{Sussman:12} (\citeyear{Sussman:12}) and
\citet{siva.spectral:11}. 

In this paper, we derive new error bounds for spectral clustering for
the purpose
of
community recovery in moderately
sparse stochastic block models and degree corrected stochastic block models
[see, e.g., \citet{KarrerN11}], where the maximum expected node degree
is of order
$\log n$ or higher. Our main contribution is to show that the most basic
form of spectral clustering is successful in recovering the latent community
memberships under conditions on the network sparsity that are weaker
than the
ones used in
most of literature.
Our results yield some sharpening of existing analyses of spectral
clustering for community recovery, and provide a theoretical
justification for the effectiveness
of this procedure in moderately sparse networks.
We take note that there are competing methods yielding consistent community
recovery under even milder conditions on the rate of growth of the node degrees,
but they either rely on
combinatorial methods that are
computationally demanding [\citet{BickelC09}] or are guaranteed to be
successful
provided that they are given good starting points [\citet{AminiCBL12}],
which are typically unknown. Other computationally efficient
procedures with strong theoretical guarantees, which include in
particular the ones
proposed and analyzed in \citet{McSherry01,ChenSX12,stephane,sarkar.bickel:13},
require instead the degrees to be of larger order than $\log n$.
More detailed comparisons with some of these contributions will
be given after the statement of main results as more technical
background is introduced.
Finally, it is also known that in the ultra-sparse case, where the
maximum degree
is of order $O(1)$, consistent community recovery is impossible and one
can only hope to recover the communities
up to a constant fraction
[see \citet
{Coja-Oghlan10,Decelle11,Krzakala13,Massoulie13}, \citeauthor{MosselNS12}
(\citeyear{MosselNS12,MosselNS13})].


The contributions of this paper are as follows. We prove that a
simplest form of
spectral clustering, consisting of applying approximate $k$-means algorithms
to the rows of
the matrix formed by the leading eigenvectors of the adjacency matrix, allows
to recover the community membeships of all but a vanishing fraction of the
nodes in stochastic block models with expected degree as small as $\log n$,
with high probability.
We also extend this result to degree corrected
stochastic block models by analyzing an approximate spherical
$k$-median spectral clustering algorithm.
The algorithms we consider are among the most practical and
computationally affordable procedures
available. Yet the theoretical guarantees we provide hold under rather general
assumptions of sparsity that are weaker than the ones used in
algorithms of
similar complexity.
Our arguments extend those in \citet{RoheCY11} and \citet{Jin12} by
combining a principal subspace perturbation analysis (Lemma~\ref
{lem:variational-bound}), a deterministic performance guarantee of
approximate $k$-means clustering (Lemma~\ref{lem:hamming}) and
a sharp bound on the spectrum of binary random matrices (Theorem~\ref
{thm:long.thm}), which
may be of independent interest.
These techniques give sharper results under weaker conditions.
In particular, the subspace perturbation analysis
allows us to avoid the individual eigengap
condition. On the other hand, the spectral bound gives a better large deviation
result that cannot be obtained by the matrix Bernstein inequality
[\citet{ChungR11,Tropp12}] and leads to a simple extension to the
degree corrected stochastic block model.


The article is organized as follows. In Section~\ref{sec:prelim} we give
formal introduction to the stochastic block model
and spectral clustering.
The main results are presented and compared to related works in
Section~\ref{sec:sbm} for regular SBM's
and in Section~\ref{sec:dcbm} for degree corrected block models.
Section~\ref{sec:general} presents the proofs of main results, including
a general, highly modular scheme of
analyzing performance of spectral clustering algorithms.
Concluding remarks are given in Section~\ref{sec:disc}.

\textit{Notation}.
For a matrix $M$ and index sets $\mathcal I, \mathcal J\subseteq[n]$,
let $M_{\mathcal I *}$ and $M_{* \mathcal J}$ be the submatrix of $M$
consisting the corresponding rows and columns.
Let $\mathbb M_{n,K}$ be the collection of all $n\times K$ matrices
where each row
has exactly one 1 and $(K-1)$ 0's. For any $\Theta\in\mathbb M_{n,K}$,
we call $\Theta$ a
\emph{membership matrix}, and
the community membership of a node $i$ is denoted by $g_i\in\{1,\ldots,K\}
$, which
satisfies
$\Theta_{i g_i}=1$.
Let $G_k=G_k(\Theta)
=\{1\le i\le n\dvtx g_i=k\}$ and $n_k=|G_k|$ for all $1\le k\le K$. Let
$n_{\min}=\min_{1\le k\le K}
n_k$, $n_{\max}=\max_{1\le k\le K} n_k$, and $n_{\max}'$ be the second
largest community size.
We use $\|\cdot\|$ to denote both the Euclidean norm of a vector and
the spectral
norm of a matrix. $\|M\|_F=(\operatorname{ trace}(M^T M))^{1/2}$ denotes the
Frobenius norm of
a matrix $M$. The $\ell_0$ norm $\|M\|_0$
simply counts the number of nonzero entries in $M$. For any square
matrix $M$, $\operatorname{ diag}(M)$
denotes the matrix obtained by setting all off-diagonal entries of $M$
to 0.
For two sequences of real numbers $\{ x_n \}$ and $\{ y_n \}$, we will
write $x_n =
o(y_n)$ if $\lim_n x_n/y_n = 0$, $x_n = O(y_n)$ if $ \llvert  x_n/y_n
\rrvert  \leq
C$ for all $n$ and some positive $C$ and $x_n = \Omega(y_n)$ if $
\llvert  x_n/y_n
\rrvert  >C$ for
all $n$ and some positive~$C$.

\section{Preliminaries}\label{sec:prelim}
\subsection{Model setup}\label{sec:setup}
A stochastic block model with $n$ nodes and $K$ communities is
parameterized by a pair of matrices
$(\Theta,B)$, where $\Theta\in\mathbb M_{n,K}$ is the membership
matrix and
$B\in\mathbb R^{K\times K}$ is a symmetric \emph{connectivity matrix}.
For each node $i$, let $g_i$ ($1\le g_i\le K$) be its community label,
such that the $i$th row of
$\Theta$ is $1$ in column $g_i$ and 0 elsewhere.
On the other hand, the entry
$B_{k\ell}$ in $B$ is the edge probability between any node in
community $k$ and any node in community $\ell$. Given $(\Theta,B)$, the
adjacency matrix $A=(a_{ij})_{1\le i,j\le n}$ is generated as
\[
a_{ij}=\cases{ %
\mbox{independent $\operatorname{
Bernoulli}(B_{g_ig_j})$},& \quad $\mbox{if $i<j$},$
\vspace*{2pt}\cr1
0,&\quad $\mbox{if $i=j$},$
\vspace*{2pt}\cr
a_{ji}, &\quad $\mbox{if $i>j$}.$}
\]
The goal of community recovery is to recover the membership matrix
$\Theta$
up to column permutations. Throughout this article, we assume that
the number of communities, $K$, is known.
For an estimate $\hat\Theta\in\mathbb M_{n,K}$ of the node memberships,
we consider two measures of estimation error.
The first one is an overall relative
error
\[
L(\hat\Theta,\Theta)=n^{-1}\min_{J\in E_{K}}\|\hat\Theta J-
\Theta \|_0,
\]
where $E_K$ is the set of all $K\times K$ permutation matrices. Because
both $\hat\Theta J$ and $\Theta$
are membership matrices, we have $\|\hat\Theta J-\Theta\|_0=\|\hat
\Theta J-\Theta\|_F^2$.
This quantity measures the overall
proportion of mis-clustered nodes.

The other performance criterion measures the worst case relative error
over all communities:
\[
\tilde L(\hat\Theta,\Theta)=\min_{J\in E_K}\max
_{1\le k\le
K}n_{k}^{-1}\bigl\|(\hat\Theta
J)_{G_k*} -\Theta_{G_k*}\bigr\|_0.
\]
%
It is obvious that $0\le L(\hat\Theta,\Theta)\le\tilde L(\hat
\Theta,\Theta)\le2$.
Thus, $\tilde L$ is a stronger criterion than $L$ in that it requires the
estimator to do well for all communities, while an estimator $\hat
\Theta
$ with
small $L(\hat\Theta,\Theta)$ may have large relative errors for some
small communities.

\subsection{Spectral clustering}
Spectral clustering is a simple method for community recovery [\citet
{VonLuxburg07,RoheCY11,Jin12}].
In a SBM, the heuristic of spectral clustering
is to relate the eigenvectors of
$A$ to those of $P:=\Theta B\Theta^T$ using the fact that
$\mathbb E(A)=P-\operatorname{ diag}(P)$.
Let $P=UDU^T$ be the eigen-decomposition of $P$ with
$U^T U=I_K$ and $D\in\mathbb R^{K\times K}$ diagonal, then
it is easy to see that $U$ has only $K$ distinct rows since $P$ has
only $K$
distinct rows.
Under mild conditions, it is also the case that
two nodes are in the same community if and only if their corresponding
rows in $U$ are
the same. This is formally stated in the following lemma.

\begin{lemma}[(Basic eigen-structure of SBMs)]\label{lem:p-eigen}
Let the pair $(\Theta,B)$ parame\-trize a SBM with $K$ communities, where
$B$ is full rank.
Let $UDU^T$ be the eigen-decomposition of $P=\Theta B\Theta^T$.
Then
$U=\Theta X$ where $X\in\reals^{K\times K}$ and
$\|X_{k*}-X_{\ell*}\|=\sqrt{n_k^{-1}+n_\ell^{-1}}$ for all
$1\le k< \ell\le K$.
\end{lemma}

\begin{pf}
Let $\Delta=\operatorname{ diag}(\sqrt{n_1},\ldots,\sqrt{n_K})$ then
%
\begin{equation}
\label{eq:eigen-p} P=\Theta B\Theta=\Theta\Delta^{-1}\Delta B \Delta\bigl(
\Theta\Delta ^{-1}\bigr)^T.
\end{equation}
It is straightforward to verify that $\Theta\Delta^{-1}$ is orthonormal.
Let $ZDZ^T= \Delta B \Delta$ be the eigen-decomposition of $\Delta B
\Delta$.
Thus, we have $P=UDU^T$ where $U=\Theta\Delta^{-1}Z$. The claim follows
by letting $X=\Delta^{-1} Z$ and realizing that the rows of $\Delta
^{-1} Z$ are perpendicular to each other
and the $k$th row has length $\|(\Delta Z)_{k*}\|=\sqrt{1/n_k}$.
\end{pf}

Based on this observation, spectral clustering tries to estimate $U$
and its row clustering
using a spectral decomposition of $A$. The intuition for the procedure
is as
follows.
Consider the difference between $A$ and $P$:
\[
A-P=\bigl(A-\mathbb E(A)\bigr)-\operatorname{ diag}(P),
\]
which is
a symmetric noise matrix plus a diagonal matrix. Intuitively, the
eigenvectors of $A$ will be close to those of $P$ because the
eigenvalues of $P$ scales linearly with $n$ while
the noise matrix $(A-\mathbb E(A))$ has operator norm on the scale of
$\sqrt{n}$ and $\operatorname{ diag}(P)$ is like a constant. Therefore, letting
$A=\hat U \hat D\hat U^T$ be the $K$-dimensional eigen-decomposition of
$A$ corresponding to the $K$ largest absolute eigenvalues, we can see
that $\hat
U$ should have roughly $K$ distinct rows because they are slightly
perturbed versions of the rows in $U$. Then one should be able to obtain
a good community partition by applying a clustering algorithm
on the rows of $\hat U$.
In this paper we consider the $k$-means clustering, defined as
%
\begin{equation}
\label{eq:k-means} (\hat\Theta,\hat X)=\arg\min_{\Theta\in\mathbb M_{n,K},X\in
\mathbb
R^{K\times K}} \|\Theta X -
\hat U\|_F^2.
\end{equation}
It is known that
finding a global minimizer for the $k$-means problem (\ref
{eq:k-means}) is
NP-hard [see, e.g., \citet{AloiseDHP09}]. However, efficient
algorithms exist for finding an approximate solution whose value is
within a constant
fraction of the optimal value [\citet{KumarSS04}]. That is, there are
polynomial time algorithms
that find
%
\begin{eqnarray}
\label{eq:k-means-eps} %
(\hat\Theta,\hat X)\in\mathbb
M_{n,K} \times\reals^{K\times K}
\nonumber
\\[-8pt]
\\[-8pt]
\eqntext{\mathrm{s.t.} \quad\displaystyle \|\hat\Theta\hat X-\hat U\|_F^2\le(1+
\varepsilon) \min_{\Theta\in
\mathbb
M_{n,K},X\in\mathbb R_{K\times K}} \|\Theta X - \hat U\|_F^2.}
\end{eqnarray}
The spectral clustering algorithm we consider here is summarized in
Algorithm \hyperref[alg1]{1}.

\begin{figure}
\begin{center}
\fbox{\parbox{4.5in}{
\begin{center}{\textsf{Algorithm 1:
Spectral clustering with approximate $k$-means}}\label{alg1}
\end{center}

\textbf{Input:} Adjacency matrix $A$; number of communities $K$;
approximation parameter $\varepsilon$.\\
\textbf{Output:} Membership matrix $\hat\Theta\in\mathbb M_{n, K}$.
\begin{enumerate}
\item Calculate $\hat U\in\reals^{n\times K}$ consisting of the leading
$k$ eigenvectors (ordered
in absolute eigenvalue) of $A$.
\item Let $(\hat\Theta,\hat X)$ be an $(1+\varepsilon)$-approximate
solution to
the $k$-means problem \eqref{eq:k-means-eps} with $K$ clusters and
input matrix $\hat U$.
\item Output $\hat\Theta$.
\end{enumerate}
}}
\end{center}
\end{figure}

\subsection{Sparsity scaling}
Real-world large scale
networks are usually sparse, in the sense that the number of edges from
a node (the node degree) are very small compared to the total number of nodes.
Generally speaking, community recovery is
hard when data is sparse. 
As a result, an important criterion of evaluating a community recovery
method is its performance under different levels of sparsity (usually measured
in the error rate as a function of the average/maximum degree).
The following prototypical example
exemplifies well the roles played by network sparsity as well as
other model parameters in determining the hardness of community recovery.

\begin{example}\label{exa:sbm}
Consider a SBM with $K$ communities parameterized by $(\Theta,B)$
where 
%
\begin{equation}
\label{eq:example_sbm} B=\alpha_n B_0;\qquad B_0=\lambda
I_K+(1-\lambda)\mathbf1_K \mathbf1_K^T,\qquad
0<\lambda<1,
\end{equation}
$I_K$ is the $K\times K$ identity matrix, and $\mathbf1_K$ is the
$K\times1$ vector of 1's.
\end{example}

Example~\ref{exa:sbm} assumes that the edge probability between any
pair of nodes
depends only on whether they belong to the same community.
In particular, the edge probability is $\alpha_n$ within community and
$\alpha_n(1-\lambda)$ between
community. The quantity $\lambda$ reflects the relative difference in
connectivity
between communities and within communities.
The network sparsity is captured by $\alpha_n$, where
$n\alpha_n$ provides an upper bound on the average (and maximum in
this example)
expected node degree. It can be easily seen that if $\alpha_n$ or
$\lambda$ are close to 0 then
it is hard to identify communities.

The hardness of community reconstruction
also depends on the number of communities
and the community size imbalance. For example, the famous planted
clique problem
concerns community recovery under a SBM with $K=2$ and
%
\begin{eqnarray}
\label{eq:pcp} B=\pmatrix{ 1 &1/2
\vspace*{2pt}\cr
1/2 & 1/2 }.
\end{eqnarray}
In the planted clique problem, it is known that community recovery is
easy if $n_1\ge c\sqrt{n}$ for a constant $c$ [see \citet{DeshpandeM13}
and references therein]
and on the other hand
no polynomial time algorithms have been found to succeed when
$n_1=o(\sqrt{n})$.

 \begin{remark*}
The primary concern of this paper is the effect of
$\alpha_n$ on the performance of spectral clustering. Nevertheless,
our results
explicitly keep track of other quantities such as $K$, $\lambda$,
$n_{\max}$ and
$n_{\min}$, all of which are allowed to change with $n$ in a
nontrivial manner.
The dependence of recovery error bound on some of these quantities,
such as $K$ and $\lambda$, is concerned by some authors, such as
\citet
{ChenSX12,ChaudhuriCT12,Anandkumar13}.
For ease of readability, we do not always make this dependence on~$n$
explicit in our notation.
\end{remark*}
\section{Stochastic block models}\label{sec:sbm}
Our main result provides an upper bound on relative community
reconstruction error
of spectral clustering for a SBM $(\Theta, B)$ in terms of several
model parameters.

\begin{theorem}\label{thm:main}
Let $A$ be an adjacency matrix generated from a stochastic block model
$(\Theta,B)$. Assume that
$P=\Theta B\Theta^T$ is of rank $K$, with smallest absolute nonzero
eigenvalue at least $\gamma_n$
and $\max_{k,\ell}B_{k\ell}\le\alpha_n$
for some $\alpha_n\ge\log n/n$.
Let $\hat\Theta$ be
the output of spectral clustering using $(1+\varepsilon)$-approximate
$k$-means (Algorithm~\hyperref[alg1]{1}).
There exists an absolute constant $c>0$, such that, if
%
\begin{equation}
\label{eq:condition.main.thm} (2 + \varepsilon)\frac{ K n \alpha_n}{\gamma_n^2} < c,
\end{equation}
then,
with probability at least $1-n^{-1}$, there exist subsets
$S_k\subset G_k$ for $k=1,\ldots,K$, and a $K\times K$ permutation matrix
$J$ such that
$\hat\Theta_{G*}J=\Theta_{G*}$, where $G=\bigcup_{k=1}^K (G_k\setminus S_k)$,
and
%
\begin{equation}
\label{eq:result-main-sbm} \sum_{k=1}^n
\frac{|S_k|}{n_k}\le c^{-1}(2+\varepsilon)\frac
{Kn\alpha
_n}{\gamma_n^2}.
\end{equation}
\end{theorem}

The proof of Theorem~\ref{thm:main}, given in Section~\ref{sec:general},
is modular, and can be derived from
several relatively independent lemmas.

The sets $S_k$ ($1\le k\le K$) consist of nodes in $G_k$ for which the
clustering correctness
cannot be guaranteed.
The permutation matrix $J$ in the above theorem leads
to an upper bound on reconstruction error $\tilde L(\hat\Theta,\Theta)$
[and hence on $L(\hat\Theta,\Theta)$] through equation (\ref
{eq:result-main-sbm}).

Condition \eqref{eq:condition.main.thm} specifies the range of
model parameters $(K,n,\gamma_n,\alpha_n)$ for which the result is applicable.
It is included only for technical reasons, because
it holds whenever the bound in \eqref{eq:result-main-sbm} vanishes and,
therefore, implies consistency.
In particular, as discussed after Corollary~\ref{cor:main}, we have
$Kn\alpha
_n/\gamma_n^2=o(1)$ in many
interesting cases.
The constant $c$ in \eqref{eq:condition.main.thm} can be written as
$c=1/(64C^2)$ where $C$ is an absolute constant defined in Theorem~\ref
{thm:long.thm} and
can be explicitly tracked
in the proof presented in the supplementary material [\citet{SUPP}]. The
assumption
of $\alpha_n \ge\log n/n$ can be changed to $\alpha_n\ge c_0\log n/n$
for any $c_0>0$, and also
the probability bound $1-n^{-1}$ can be changed to $1-n^{-r}$ for any $r>0$,
with a different constant $c=c(c_0,r)$ in \eqref{eq:condition.main.thm}
and \eqref{eq:result-main-sbm}.

While Theorem~\ref{thm:main} provides a general error bound for
spectral clustering,
the quantities involved are not in the most transparent form.
For example, the bound does not clearly reflect the intuition that the
error should increase when $\alpha_n$ decreases.
This is because the quantity $\gamma_n$ contains the parameter $\alpha
_n$. Also the dependence on the community
size imbalance as well as the community separation (which corresponds to
the parameter $\lambda$ in Example~\ref{exa:sbm}) remains unclear.
The next corollary illustrates the error bound in terms of these model
parameters.

\begin{corollary}\label{cor:main}
Let $A$ be an adjacency matrix from the SBM $(\Theta,B)$, where
$B=\alpha_nB_0$ for some $\alpha_n\ge\log n/n$ and with $B_0$ having
minimum absolute
eigenvalue $\ge\lambda>0$ and $\max_{k\ell}B_0(k,\ell)=1$.
Let $\hat\Theta$ be
the output of spectral clustering using $(1+\varepsilon)$-approximate
$k$-means (Algorithm~\hyperref[alg1]{1}). Then there exists an absolute constant $c$
such that
if
%
\begin{equation}
\label{eq:condition-cor-main} (2+\varepsilon)\frac{K n}{n_{\min}^2
\lambda^2 \alpha_n}< c
\end{equation}
then with probability at least $1-n^{-1}$,
\[
\tilde L(\hat\Theta,\Theta)\le c^{-1}(2+\varepsilon)\frac
{Kn}{n_{\min
}^2\lambda^2\alpha_n}
\]
and
\[
L(\hat\Theta,\Theta)\le c^{-1}(2+\varepsilon)\frac{Kn_{\max
}'}{n_{\min
}^2\lambda^2\alpha_n}.
\]
\end{corollary}

In the special case of a balanced community sizes
[i.e., $n_{\max}/n_{\min}=O(1)$] and constant $\lambda$,
if $\alpha_n=\Omega(\log n/n)$, then $L(\hat\Theta,\Theta)
=O_P(K^2(n\alpha_n)^{-1})=O_P(K^2/ \log n)$. Thus $L(\hat\Theta,\Theta)
=o_P(1)$ if $K=o(\sqrt{\log n})$. This improves the results in \citet
{RoheCY11}
where $\alpha_n$ needs to be of order $1/\log n$ for a similar result.

In Example~\ref{exa:sbm}, the smallest nonzero eigenvalue of $B_0$ is
$\lambda
$. Recall that
$\lambda$ is the relative difference of within- and between-community
edge probabilities. Corollary~\ref{cor:main} then implies that when this
relative difference stays bounded away from zero,
the communities can be consistently recovered by simple
spectral clustering
as long as the expected node degrees are no less than $\log n$. On the
other hand,
when $\alpha_n$ is constant and $\lambda=\lambda_n$ varies with $n$,
spectral clustering can recover the
communities when the relative edge probability gap grows faster than
$1/\sqrt{n}$.

In the planted clique problem, $L(\hat\Theta,\Theta)$ has limited
meaning because
a trivial clustering putting all nodes in one cluster achieves $L(\hat
\Theta,\Theta)=2n_{\min}/n$
which is $o(1)$ in the most interesting regime.\vspace*{1pt} Therefore, it makes
more sense to
consider $\tilde L(\hat\Theta,\Theta)$.
Now $B_0=B$ is given by \eqref{eq:pcp}, with minimum eigenvalue $>0.19$.
Applying Corollary~\ref{cor:main} with
$K=2$, $\lambda=0.19$, $\alpha_n=1$, and any fixed $\varepsilon>0$,
we have
\[
\tilde L(\hat\Theta,\Theta)< c' \frac{n}{n_{\min}^2},
\]
provided that $c'n/n_{\min}^2<1$, where $c'$ is a different absolute constant.
Therefore, when $n_{\min}\ge\sqrt{an}$ for some $a>c'$,
$\hat\Theta$ recovers the hidden clique with a relative error no larger
than $c'/a$.
Thus, our result reaches the well believed computation barrier [up to
constant factor, see
\citet{DeshpandeM13} and references therein]
of the planted clique problem. 

There are spectral methods other than
spectral clustering that can provide consistent community recovery. One such
well-known
method is the procedure analyzed by \citet{McSherry01}.
The planted partition
problem in that setting corresponds to the problem of recovering the community
memberships in the
SBM. To simplify the presentation and focus on the dependence of
network sparsity, we
consider the SBM in Example~\ref{exa:sbm} with two equal-sized
communities and
a constant $\lambda\in(0,1)$. According to
Theorem~4 in \citet{McSherry01}, that method can recover
the true communities with probability at least $ 1- n^{-1}$ provided
that, after some simplification,
%
\begin{equation}
\label{eq:mcsherry} \lambda^2 \alpha_n^2 n > c
\sigma_n^2 \log n \quad \mbox{and}\quad \sigma_n^2>(
\log n)^6 / n ,
\end{equation}
for some constant $c$, where
$\sigma_n^2$ is an upper bound on the maximal variance of the edges.
Therefore, the condition \eqref{eq:mcsherry}
implies that
$\alpha_n > \sqrt{c}\lambda^{-1} (\log n)^{3.5} / n$, which is
stronger than the condition in our Corollary~\ref{cor:main}.




\section{Degree corrected stochastic block models}\label{sec:dcbm}
The degree corrected
block model [DCBM, \citet{KarrerN11}] extends the standard SBM by introducing
node specific parameters to allow for varying degrees even within the same
community. A DCBM is parameterized by a triplet $(\Theta,B,\psi)$, where,
in addition to the membership matrix $\Theta$ and connectivity matrix $B$,
the vector $\psi\in\mathbb R^n$ is included to model additional
variability of
the edge probabilities at the node level. Given $(\Theta,B,\psi)$, the
edge probability
between nodes $i$ and $j$ is $\psi_i\psi_j B_{g_i g_j}$ (recall that
$g_i$ is the community label of node $i$). Similar to the SBM, the DCBM
also assumes independent edge formation given $(\Theta, B,\psi)$.
The inclusion of $\psi$ raises an issue of identifiability. So we assume
that $\max_{i\in G_k}\psi_i=1$ for all $k=1,\ldots,K$. The SBM can
be viewed as a special case of DCBM with $\psi_i=1$ for all $i$. The
DCBM greatly enhances the
flexibility of modeling degree heterogeneity and is
able to fit network data with arbitrary degree distribution. Successful
application
and theoretical developments
can be found in \citet{ZhaoLZ12} for likelihood methods, and in \citet
{ChaudhuriCT12,Jin12}
for spectral methods.
%

\textit{Additional notation about the degree heterogeneity}.
Let $\phi_k$ be the $n\times1$ vector that agrees with $\psi$ on $G_k$
and zero otherwise. Define $\tilde\phi_k=\phi_k/\|\phi_k\|$ and
$\tilde\psi=\sum_{k=1}^K \tilde\phi_k$.
Let $\tilde\Theta$ be a normalized membership matrix such that
$\tilde\Theta(i,k)=\tilde\psi_i$ if $i\in G_k$ and $\tilde\Theta
(i,k)=0$ otherwise.
We also define
\emph{effective community size}
$\tilde n_k:=\|\phi_k\|^2$.
Let $\tilde n_{\min}=\min_k \tilde n_k$
and $\tilde n_{\max}=\max_k \tilde n_k$.

The spectral clustering heuristic can be extended to DCBMs
by considering the eigen-decomposition $P=UDU^T$
where $P=\operatorname{ diag}(\psi)\Theta B \Theta^T\operatorname{ diag}(\psi)$.
Now the matrix $U$ may
have more than $K$ distinct rows due to the effect of $\psi$. However,
the rows of $U$ point to at most $K$ distinct directions [\citet{Jin12}].
The following lemma is the analogue of Lemma~\ref{lem:p-eigen} for DCBMs.

\begin{lemma}[(Spectral structure of mean matrix in DCBM)]\label
{lem:dcmb-eigen-row-norm}
Let $UDU^T$ be the eigen-decomposition of $P=\operatorname{ diag}(\psi)\Theta B
\Theta^T\operatorname{ diag}(\psi)$ in a DCBM parameterized
by $(\Theta,B,\psi)$. Then there exists a $K\times K$ orthogonal matrix
$H$ such that
\[
U_{i*}=\tilde\psi_{i}H_{k*}\qquad\forall 1\le k\le K,
i\in G_k.
\]
\end{lemma}

\begin{pf}
First, realize that $\operatorname{ diag}(\psi)\Theta=\tilde\Theta\Psi$,
where $\Psi=\operatorname{ diag}(\|\phi_1\|,\ldots,\break  \|\phi_K\|)$.
%
\begin{equation}
\label{eq:eigen-P-dcbm} P=\operatorname{ diag}(\psi)
\Theta B\Theta^T\operatorname{ diag}(\psi)=
\tilde\Theta \Psi B\Psi\tilde\Theta^T = \tilde\Theta H D (\tilde\Theta
H)^T,
\end{equation}
where $\Psi B\Psi=H D H^T$ is the eigen-decomposition of $\Psi B\Psi$.
Note that $\tilde\Theta^T\tilde\Theta= I_K$ so
$\tilde\Theta H D (\tilde\Theta H)^T$ is an eigen-decomposition of $P$.
\end{pf}

As a result, finding the true community partition corresponds to clustering
the \emph{directions} of the row vectors in $U$, where
some form of normalization must be employed in order to filter out
the nuisance parameter $\psi$.
In particular, we consider \emph{spherical clustering}, which looks for
a cluster structure
among the rows of a normalized matrix $U'$ with $U'_{i*}=U_{i*}/\|
U_{i*}\|$. 

In addition to the overall sparsity, the difficulty of community
recovery in a DCBM is also affected by small entries
of $\psi$.
Intuitively, if $\psi_i\approx0$, then it is hard to identify
the community membership of node $i$ because few edges are observed for
this node.
However, the interaction between small entries of $\psi$ and the
overall network sparsity
(the maximum/average degree)
has not been well understood.
In the analysis of profile likelihood methods,
\citet{ZhaoLZ12} assume that the entries of $\psi$ are fixed constants.
In spectral clustering, \citet{Jin12} allows milder conditions on
$\psi$ but
needs the average degree to be polynomial in $n$.

Our analysis uses the following quantity as a
summarizing measure of node heterogeneity in each community $G_k$:
\[
\nu_k:= n_k^{-2}\sum
_{i\in G_k}\tilde\psi_{i}^{-2} ,\qquad k=1,2,
\ldots,K.
\]
By definition $\nu_k\in[1,\infty)$ and a larger $\nu_k$ indicates %
a stronger heterogeneity
in the $k$th community.
On the other hand,
$\nu_k=1$ indicates within-community homogeneity ($\psi_i=1$ for all
$i\in G_k$).

The argument developed for SBMs in previous sections can be extended to cover
very general
degree corrected models. In particular, let $\hat U\in\reals^{n\times
K}$ consist
the $K$ leading eigenvectors of $A$.
We consider the following
spherical $k$-median spectral clustering:
%
\begin{equation}
\mathrm{ minimize}_{\Theta\in\mathbb M_{n,K}, X\in\mathbb R^{K\times K}}
\bigl\|\Theta X-\hat U'\bigr\|_{2,1}
,\label{eq:k-median}
\end{equation}
where $\hat U'$ is the row-normalized version of $\hat U$
and
$\|M\|_{2,1}=\sum_{i=1}\|M_{i*}\|$ is the matrix $(2,1)$-norm.
We will not
require to solve \eqref{eq:k-median} exactly but instead we consider a
$(1+\varepsilon)$ approximation
$(\hat\Theta,\hat X)$ to the $k$-median
problem, which can be solved in polynomial time when $\varepsilon>
\sqrt{3}$
[\citet{CharikarGTS99,LiS13}]. The practical procedure will also take
care of the possible zero rows in
$\hat U$ and is described in detail in Algorithm \hyperref[alg2]{2}.

\begin{figure}\label{alg2}
\begin{center}
\fbox{\parbox{4.5in}{
\begin{center}{\textsf{Algorithm 2:
Spherical $k$-median spectral clustering}}
\end{center}

\textbf{Input:} Adjacency matrix $A$; number of communities $K$;
approximation parameter $\varepsilon$.\\
\textbf{Output:} Membership matrix $\hat\Theta\in\mathbb M_{n,K}$.
\begin{enumerate}
\item Calculate $\hat U\in\reals^{n\times K}$ consisting of the leading
$k$ eigenvectors (ordered
in absolute eigenvalue) of $A$.
\item Let $I_+=\{i\dvtx\|\hat U_{i*}\|>0\}$ and $\hat U^{+}=(\hat U_{I_+*})$.
\item Let $\hat U'$ be row-normalized version of $\hat U^+$.
\item Let $(\hat\Theta^+,\hat X)$ be an $(1+\varepsilon)$-approximate
solution to
the $k$-median problem with $K$ clusters and input matrix $\hat U'$.
\item Output $\hat\Theta$ with $\hat\Theta_{i*}$ being the
corresponding row in
$\hat\Theta^+$ if $i\in I_+$, and $\hat\Theta_{i*}=(1,0,\ldots,0)$ if
$i\notin I_+$.
\end{enumerate}
}}
\end{center}
\end{figure}

\subsection{Analysis of spherical $k$-median spectral clustering for DCBM}
%
We have the following main theorem for spherical $k$-median spectral
clustering in
DCBMs. It is proved in Appendix~\ref{sec:pf-dcbm}.

\begin{theorem}[(Main result for DCBM)]\label{thm:main-dcbm}
Consider a DCBM $(\Theta,B,\psi)$ with $K$ communities, where $P=\operatorname{
diag}(\psi) \Theta B\Theta^T \operatorname{ diag}(\psi)$
has rank $K$, the smallest nonzero absolute eigenvalue at least
$\gamma_n$, and the maximum entry
bounded from above by $\alpha_n\ge\log n/n$.
There exists an absolute constant $c > 0$ such that if
%
\begin{equation}
\label{eq:dcbm-main-condition} (2.5+\varepsilon)
\frac{\sqrt{Kn\alpha_n}}{\gamma_n}< c\frac{n_{\min}}{\sqrt{\sum_{k=1}^K n_k^2 \nu_k}}
\end{equation}
then,
with probability at least $1-n^{-1}$,
%
\begin{equation}
\label{eq:dcbm-main-rate} L(\hat\Theta,\Theta)\le c^{-1}(2.5+\varepsilon)\sqrt
{\sum_{k=1}^K n_k^2
\nu _k}\frac{\sqrt{K\alpha_n}}{\gamma_n\sqrt{n}}.
\end{equation}
\end{theorem}

\begin{remark*} The constant $c$ equals $1/(8C)$ where $C$ is
the universal constant
in Theorem~\ref{thm:long.thm}. The condition on $\alpha_n$
and probability guarantee can also be changed to $\alpha_0\ge c_0\log
n/n$ and $1-n^{-r}$,
respectively, with a different constant $c=c(c_0,r)$ in equations (\ref
{eq:dcbm-main-condition}) and (\ref{eq:dcbm-main-rate}).
\end{remark*}

Theorem~\ref{thm:main-dcbm} immediately implies a counterpart of
Corollary~\ref{cor:main} under more explicit scaling of the model parameters.

\begin{corollary}\label{cor:main-dcbm}
Let $A$ be an adjacency matrix from DCBM $(\Theta,B,\psi)$, such that
$B=\alpha_nB_0$ for some $\alpha_n\ge\log n/n$ where $B_0$ has minimum
absolute eigenvalue $\lambda>0$ and $\max_{k\ell}B_0(k,\ell)=1$. Let
$(\hat\Theta,\hat X)$ be an $(1+\varepsilon)$-approximate solution
to the spherical $k$-median algorithm (Algorithm \hyperref[alg2]{2}).
There exists an absolute constant $c$ such that if
\[
(2.5+\varepsilon)\frac{\sqrt{Kn}}{\tilde n_{\min} \lambda
\sqrt{\alpha_n}}< c\frac{n_{\min}}{\sqrt{\sum_{k=1}^K n_k^2 \nu
_k}} ,
\]
then, with probability at least $1-n^{-1}$,
\[
L(\hat\Theta,\Theta)\le c^{-1}(2.5+\varepsilon)\frac{\sqrt
{K}}{\tilde
n_{\min}\lambda\sqrt{n\alpha_n}}\sqrt
{\sum_{k=1}^K n_k^2
\nu_k}.
\]
\end{corollary}

Comparing with Theorem~\ref{thm:main} and Corollary~\ref{cor:main},
the results for DCBM
are different in two major aspects. First, the DCBM condition \eqref
{eq:dcbm-main-condition} involves the term
$n_{\min}^2/\sum_{k=1}^K n_k^2\nu_k$ which is smaller than 1 (indeed
smaller than
$1/K$). This makes \eqref{eq:dcbm-main-condition} more stringent than
\eqref{eq:condition.main.thm}. Also the upper bound on $L(\hat\Theta
,\Theta)$ is
different in the same manner. Furthermore,
the argument used to prove Theorem~\ref{thm:main-dcbm} is not likely to
provide a sharp upper bound on $\tilde L(\hat\Theta,\Theta)$. We
believe this has
to do with the additional normalization step used in the spherical
$k$-median algorithm
as well as the specific strategy used in our proof.

To better understand this result, consider Example~\ref{exa:sbm} with balanced
community size: $n_{\max}/n_{\min}=O(1)$. To work with a DCBM, assume
in addition that the node degree vector $\psi$
has comparable degree heterogeneity across communities: $c_1 \nu\le
\nu
_k\le c_2 \nu$ for constants $c_1$, $c_2$. Then Corollary~\ref{cor:main-dcbm}
implies an overall relative error rate
%
\begin{equation}
L(\hat\Theta,\Theta)=O_P \biggl( \frac{\sqrt{\nu}}{\tilde n_{\min}\lambda\sqrt{n\alpha_n}} \biggr)
.\label{eq:example-rate-dcbm}
\end{equation}
Several observations are worth mentioning.
First, the error rate depends on $\nu$, the degree heterogeneity
measure, in a simple manner.
Second, the community size $n_{\min}$ that appears in Corollary~\ref{cor:main}
is replaced by $\tilde n_{\min}=\min_{k}\|\phi_k\|$, the minimum
effective sample size. Roughly speaking,
$\tilde n_{\min}\asymp n_{\min}$ as long as a constant fraction of
nodes have their $\psi_i$'s bounded away
from zero (but the rest should not be too small in order to keep $\nu$ small).
Third, if there is no degree heterogeneity ($\nu_k\equiv1$ and
$\tilde
n_{\min}=n_{\min}$),
then the rate in \eqref{eq:example-rate-dcbm} is the square root of
that given by Corollary~\ref{cor:main}.
This is due to the additional normalization step (which is not
necessary since $\nu=1$) involved in
spherical $k$-median
and the different argument used to analyze the spherical $k$-median algorithm.
Moreover, the relative error can still be $o_P(1)$ even when
$\alpha_n$ is as small as $\log n /n$, provided that $1/\nu$, $\tilde
n_{\min}/n$, and $\lambda$
stay bounded away from zero or approach zero sufficiently slowly.

\subsubsection*{Comparisons with existing work} There are relatively fewer
results for community recovery in degree corrected block models that
allow the
maximum node degree to be
of order $o(n)$.
\citet{ChaudhuriCT12} extended the method of \citet{McSherry01} to degree
corrected
block models. In the setting of Example~\ref{exa:sbm} with equal community
size, their main result (Theorems 2 and 3 in their paper)
requires $\alpha_n$ to be at least of order $1/\sqrt{n}$.
A similar requirement of a polynomial growth of expected average degree
is implicitly imposed in \citet{Jin12}, who first studied
the performance of normalized $k$-means spectral clustering in degree
corrected block models.




\section{Proof of the main results}\label{sec:general}
In this section, we present a general scheme to prove
error bounds for spectral clustering. It
contains the SBM as a special case and can be easily extended to
the degree corrected block model.
Our argument consists of three parts: (1) control the perturbation of
principal subspaces for general symmetric matrices, (2) bound the
spectrum of random binary matrices,
and (3) error bound of
$k$-mean and spherical $k$-median clustering.

\subsection{Principal subspace perturbation}
The first ingredient of our proof is to bound the difference between
the eigenvectors of
$A$ and those of $P$, where $A$ can be viewed as a noisy version of $P$.

\begin{lemma}[(Principal subspace perturbation)]
\label{lem:variational-bound}
Assume that $P\in\mathbb R^{n\times n}$ is a
rank $K$ symmetric matrix with smallest nonzero singular value $\gamma_n$.
Let $A$ be any symmetric matrix and
$\hat U, U\in\mathbb R^{n\times K}$ be the $K$ leading
eigenvectors of $A$ and $P$, respectively.
Then there exists a $K\times K$ orthogonal matrix $Q$ such that
\[
\|\hat U - U Q\|_F\le\frac{2\sqrt{2K}}{\gamma_n}\|A-P\|.
\]
\end{lemma}

Lemma~\ref{lem:variational-bound} is
proved in Appendix~\ref{sec:pf-variational}, which is based on an application
of the Davis--Kahan $\sin\Theta$ theorem [Theorem VII.3.1 of \citet
{Bhatia97}].
The presence of a $K\times K$ orthonormal matrix $Q$
in the statement of Lemma~\ref{lem:variational-bound} is to take care
of the
situation where some
leading eigenvalues have multiplicities larger than one. In this case, the
eigenvectors are determined only up to a rotation.

%
\subsection{Spectral bound of binary symmetric random matrices}
The next theorem provides a sharp probabilistic upper bound on $\|A-P\|
$ when
$A$ is a random adjacency matrix with $\E(a_{ij})=p_{ij}$.

\begin{theorem}[(Spectral bound of binary symmetric random matrices)]
\label{thm:long.thm}
Let $A$ be the adjacency matrix of a random graph on $n$ nodes in
which edges occur
independently. Set $\mathbb{E}[A] = P = (p_{ij})_{i,j=1,\ldots,n}$ and
assume that
$n\max_{ij} p_{ij}\le d$ for $d \ge c_0\log n$ and $c_0 > 0$.
Then, for any $r>0$ there exists a
constant $C = C(r,c_0)$ such that
\[
\|A-P\|\le C\sqrt{d}
\]
with probability at least $1-n^{-r}$.
\end{theorem}

This result does not follow conventional matrix concentration inequalities
such as the matrix Bernstein inequality (which will only give $\sqrt{d
\log n}$).
\citet{LuP12} use a path counting technique in random matrix theory
to prove a bound of the same order but require a maximal degree $d\ge
c_0(\log n)^4$.

The proof of Theorem~\ref{thm:long.thm} is technically involved, as it uses
combinatorial arguments in order to derive
spectral bounds for sparse random matrices. Our proof is based on
techniques developed by \citet{FeigeO05} for bounding the second largest
eigenvalue of an Erd\"{o}s--R\'{e}yni random graph with edge probability
$d/n$.
The full proof is provided in \citet{SUPP}. Here we give a brief outline
of the three major steps.

\textit{Step \textup{1:} Discretization}. We first reduce controlling $\|A-P\|$
to the
problem of bounding
the supremum of $|x^T(A-P)y|$ over all pairs of vectors $x,y$ in a
finite set of grid points. For
any given pair $(x,y)$ in the grid, the quantity $x^T(A-P)y$ is decomposed
into the sum of two parts. The first part corresponds to the small
entries of
both $x$ and $y$, called \textit{light pairs}, the other
part corresponds to the larger entries of $x$ or $y$, the \textit
{heavy pairs}.

\textit{Step \textup{2:} Bounding the light pairs}. The next step is to use Bernstein's
inequality and the union bound to control the contribution of the
light pairs,
uniformly over the points in the grid.

\textit{Step \textup{3:} Bounding the heavy pairs}. In the final step,
the contribution from the heavy pairs, which cannot be simply bounded
by conventional
Bernstein's inequality, will be bounded using a combinatorial
argument on the event that the edge numbers in a collection of
subgraphs do not
deviate much from their expectation. A sharp large deviation bound for sums
of independent Bernoulli
random variables [Corollary A.1.10 of \citet{AlonS04}] is used to achieve
better rate than standard Bernstein's inequality.

\subsection{Error bound of $k$-means/$k$-median on perturbed eigenvectors}
Spectral clustering (or spherical spectral clustering) applies a
clustering algorithm
to a matrix consisting of the eigenvectors of $A$, which is close
(in view of Lemma~\ref{lem:variational-bound}
and Theorem~\ref{thm:long.thm}) to a matrix whose rows can be
perfectly clustered.
We would like to bound the clustering error in terms of the closeness between
the real input matrix $\hat U$ and the ideal input matrix $U$.

The next lemma generalizes an argument used in \citet{Jin12} and
provides an error bound for any $(1+\varepsilon)$-approximate
$k$-means solution.

\begin{lemma}[(Approximate $k$-means error bound)]
\label{lem:hamming} For $\varepsilon>0$ and any two matrices $\hat U,
U\in
\mathbb R^{n\times K}$ such that
$U=\Theta X$ with $\Theta\in\mathbb M_{n,K}$, $X\in\mathbb
R^{K\times K}$,
let $(\hat\Theta,\hat X)$ be a $(1+\varepsilon)$-approximate
solution to
the $k$-means problem in equation (\ref{eq:k-means}) and $\bar U =
\hat\Theta
\hat X$.
For any $\delta_k\le\min_{\ell\neq k}\|X_{\ell*}-X_{k*}\|$, define
$S_k=\{i\in G_k(\Theta)\dvtx\|\bar U_{i*}-U_{i*}\|\ge\delta_{k}/2\}$
then
%
\begin{equation}
\label{eq:hamming00} \sum_{k=1}^K
|S_k|\delta_k^2\le4(4+2\varepsilon)\|\hat
U-U\|_F^2.
\end{equation}
Moreover, if
%
\begin{equation}
\label{eq:moreover} (16+8\varepsilon)\|\hat U-U\|_F^2/
\delta_k^2<n_k \qquad\mbox{for all } k,
\end{equation}
then there exists
a $K\times K$ permutation matrix $J$ such that $\hat\Theta
_{G*}=\Theta_{G*}J$,
where $G=\bigcup_{k=1}^K (G_k\setminus S_k)$.
\end{lemma}

Lemma~\ref{lem:hamming} provides a performance guarantee for approximate
$k$-means clustering
under a deterministic Frobenius norm condition on the input matrix.
As suggested by a referee, the proof of Lemma~\ref{lem:hamming} shares some
similarities with the proof of Theorem~3.1 in \citet{AwasthiS12} [see
also \citet{KumarK10}], though our assumptions are slightly different. For
completeness we provide a short and self-contained proof of Lemma~\ref
{lem:hamming}
in Appendix~\ref{sec:pf-lem-hamming}, giving explicit constant factors
in the
result. 

\subsection{Proof of main results for SBM}\label{sec:pf-main-sbm}
We first prove Theorem~\ref{thm:main}.
\begin{pf*}{Proof of Theorem~\ref{thm:main}}
Combining Lemma~\ref{lem:variational-bound} and Theorem~\ref
{thm:long.thm}, we obtain
that, for some $K$-dimensional orthogonal matrix $Q$,
%
\begin{equation}
\label{eq:step1} \|\hat U - U Q\|_F\le\frac{2\sqrt{2K}}{\gamma_n}\|A-P\| \leq
\frac{2 \sqrt{2K}}{\gamma_n} C\sqrt{n \alpha_n},
\end{equation}
with probability at
least $1 - n^{-1}$, where $C$ is the absolute constant involved in
Theorem~\ref{thm:long.thm}.
(Notice that the term $d$ in Theorem~\ref{thm:long.thm} becomes $n
\alpha_n$
in the current
setting.)

The main strategy for the rest of the proof is to apply Lemma~\ref
{lem:hamming} to
$\hat{U}$ and $UQ$. To that end, Lemma~\ref{lem:p-eigen} implies that
$UQ=\Theta X Q=\Theta X'$ where $\|X'_{k*}-X'_{\ell*}\|=
\sqrt{\frac{1}{n_k}+\frac{1}{n_\ell}}$. As a result, we can choose
$\delta_k=\sqrt{1/n_k+\frac{1}{\max\{n_\ell:\ell\neq k\}}}$ in
Lemma~\ref{lem:hamming} and hence $n_k\delta_k^2\ge1$ for all $k$.
Using \eqref{eq:step1}, a sufficient condition for
\eqref{eq:moreover} to hold is
%
\begin{equation}
\label{eq:thm.main.result} (16 + 8 \varepsilon)8 C^2 K \frac{n \alpha_n}{\gamma_n^2} \leq1
\le \min_{1\le k
\le K}n_k \delta_k^2,
\end{equation}
so that \eqref{eq:condition.main.thm} indeed implies
\eqref{eq:moreover} by setting $c = \frac{1}{64 C^2}$.
In detail, the choice of $\delta_k=1/\sqrt{n_k}$ together with \eqref
{eq:hamming00}
yields that
\[
\sum_{k=1}^K |S_k| \biggl(
\frac{1}{n_k}+\frac{1}{\max\{n_\ell:\ell
\neq k\}
} \biggr) = \sum_{k=1}^K
|S_k|\delta_k^2\le 4(4+2\varepsilon)\|\hat
U-UQ\|_F^2,
\]
which, combined with \eqref{eq:step1}, gives \eqref{eq:result-main-sbm}:
\[
\sum_{k=1}^K \frac{|S_k|}{n_k} \leq4(4+2
\varepsilon) 8 C^2 \frac{K n
\alpha_n
}{\gamma_n^2} = c^{-1} (2 +
\varepsilon) \frac{K n \alpha_n
}{\gamma_n^2}.
\]
Since Lemma~\ref{lem:hamming} ensures that the membership is correctly
recovered
outside of $\bigcup_{1\le k\le K}S_k$, the claim follows.
\end{pf*}

\begin{pf*}{Proof of Corollary~\ref{cor:main}}
It is easy to see, for example, from \eqref{eq:eigen-p},
that in this specific stochastic block model setting, $\gamma_n =
n_{\min}
\alpha_n \lambda$. Then the proof of Theorem~\ref{thm:main} applies
with $\gamma_n=n_{\min}
\alpha_n \lambda$ and gives
\[
\sum_{k=1}^K |S_k| \biggl(
\frac{1}{n_k}+\frac{1}{\max\{n_\ell:\ell
\neq k\}
} \biggr)\le64 C^2 (2+\varepsilon)
\frac{K n}{n_{\min}^2\lambda
^2\alpha_n} ,
\]
which implies that
\[
\tilde L(\hat\Theta,\Theta)\le\max_{1\le k\le K}\frac
{|S_k|}{n_k}\le
\sum_{1\le k\le K}\frac{|S_k|}{n_k}\le64C^2(2+
\varepsilon)\frac{K
n}{n_{\min}^2\lambda^2\alpha_n} ,
\]
and, recalling that $n_{\max}'$ is the second largest community size,
\[
L(\hat{\Theta},\Theta)\le\frac{1}{n} \sum_{k=1}^K
|S_k| \leq64C^2(2+\varepsilon) \frac{K n_{\max}'
}{n_{\min}^2 \lambda^2\alpha_n }.
\]
%
\upqed\end{pf*}

%

\section{Concluding remarks}\label{sec:disc}
The analysis in this paper applies directly to the eigenvectors
of the adjacency matrix, by combining tools in
subspace perturbation and spectral bounds of binary random graphs.
In the literature, spectral clustering using the graph Laplacian or its
variants is very popular and can sometimes lead to better empirical performance
[\citet{VonLuxburg07,RoheCY11,sarkar.bickel:13}].
An important future work would be to
extend some of the results and techniques in this paper
to spectral clustering using the graph Laplacian.
The graph Laplacian normalizes the adjacency matrix
by the node degree, which can introduce extra noise if the
network is sparse and many node degrees are small.
In several recent works, \citet{ChaudhuriCT12,QinR13} studied
graph Laplacian based spectral clustering with regularization,
where a small constant is added to all node degrees prior to the normalization.
Further understanding the bias-variance trade off would be both
important and interesting.

For degree corrected block models,
regularization methods may also lead to error bounds with
better dependence on small entries of $\psi$. The intuition is
that $\nu_k$ can be very large even when only one $\psi_i$ is very
close to zero.
In this case, one should be able to simply discard nodes like this and
work on those
with large enough
degrees. Finding the correct regularization to diminish the effect of
small-degree nodes
and analyzing the new algorithm will be pursued in future work.

This paper aims at understanding the performance of spectral clustering
in stochastic block models. While our main focus is the performance of
spectral clustering
as the network sparsity changes,
the resulting error bounds explicitly keep track of five independent
model parameters ($K$, $\alpha_n$, $\lambda$, $n_{\min}$, $n_{\max}$).
Existing results usually develop error bounds depending on a subset of
these parameters,
keeping others as constant [see, e.g., \citet
{BickelC09,ChenSX12,ZhaoLZ12}].
In the planted clique model, our result implies that spectral clustering
can find the hidden clique when its size is at least $c\sqrt{n}$ for
some large enough
constant $c$.
Our result also provides good insight
in understanding the impact of the number of clusters and separation between
communities. For instance, in Example~\ref{exa:sbm}, let
$\alpha_n\equiv1$, $n_{\max}=n_{\min}=n/K$. Then Corollary~\ref{cor:main}
implies that
spectral clustering is consistent if
$
K^2/(n\lambda^2)\rightarrow0 $.
More generally, the guarantees of Corollary~\ref{cor:main} compares favorably
against most existing results as summarized in \citet{ChenSX12}, in
terms of allowable cluster
size, density gap and overall sparsity.
It would be interesting to develop a unified theoretical framework
(e.g., minimax theory) such that all methods and model parameters can
be studied and compared together.

\begin{appendix}
\section*{Appendix: Technical proofs}\label{appendix-a}
For any two matrices $A$ and $B$ of the same dimension, we use the notation
$\langle A, B \rangle=\operatorname{ trace}(A^T B)$ for the standard matrix
inner product.

\subsection{Proof of Lemma \texorpdfstring{\protect\ref
{lem:variational-bound}}{5.1}}\label{sec:pf-variational}
By Proposition~2.2 of \citet{VuL13}, there exists a $K$-dimensional
orthogonal matrix
$Q$ such that
\[
\frac{1}{\sqrt{2K}}\|\hat U-U Q\|_F \le\frac{1}{\sqrt{K}}\bigl\|\bigl(I-
\hat U \hat U^T\bigr)UU^T\bigr\|_F \le\bigl\|\bigl(I-
\hat U \hat U^T\bigr)UU^T\bigr\|.
\]
Next, we establish that $\|(I-\hat U \hat U^T)UU^T\|\le2\frac{\|A-P\|
}{\gamma_n}$.
If $\|A-P\|\le\gamma_n/2$, then by Davis--Kahan $\sin\Theta$ theorem,
we have
\[
\bigl\|\bigl(I-\hat U \hat U^T\bigr)UU^T\bigr\|\le
\frac{\|A-P\|}{\gamma_n-\|A-P\|}\le 2\frac{\|A-P\|}{\gamma_n}.
\]
If $\|A-P\|>\gamma_n/2$, then
\[
\bigl\|\bigl(I-\hat U \hat U^T\bigr)UU^T\bigr\|\le1\le2
\frac{\|A-P\|}{\gamma_n}.
\]

\subsection{Proof of Lemma \texorpdfstring{\protect\ref{lem:hamming}}{5.3}}\label
{sec:pf-lem-hamming}
First, by the definition of $\bar U$ and the fact that $U$ is feasible
for problem (\ref{eq:k-means}), we have
$\|\bar U - U\|_F^2\le2\|\bar U - \hat U\|_F^2+2\|\hat U - U\|_F^2\le
(4+2\varepsilon)\|\hat U-U\|_F^2$.
Then
%
\begin{equation}
\label{eq:hamming1} \sum_{k=1}^K
|S_k|\delta_k^2/4\le\|\bar U-U
\|_F^2\le(4+2\varepsilon )\| \hat U - U
\|_F^2,
\end{equation}
which concludes the first claim of the lemma.

Under the assumption described in the second part of the lemma,
equation (\ref{eq:hamming1}) further implies that
\[
|S_k|\le(16+8\varepsilon)\|\hat U-U\|_F^2/
\delta_k^2<n_k \qquad\mbox{for all } k.
\]
Therefore, $T_k\equiv G_k\setminus S_k\neq\varnothing$, for each $k$.
If
$i\in T_k$ and $j\in T_\ell$ with $k\neq\ell$, then
$\bar U_{i*}\neq\bar U_{j*}$ because
otherwise
$\max(\delta_k,\delta_\ell)\le\|U_{i*}-U_{j*}\|\le\|U_{i*}-\bar
U_{i*}\|+\|U_{j*}-\bar U_{j*}\|<\delta_k/2 +\delta_\ell/2$, which is
impossible.
This further implies that $\bar U$ has exactly $K$ distinct rows,
because the number of distinct rows is no larger than $K$ as part of the
constraints of the optimization problem \eqref{eq:k-means}.

On the other hand, if $i$ and $j$ are both in $T_k$, for some
$k$, then $\bar U_{i *} = \bar U_{j *}$ because
otherwise there would be more than $K$ distinct rows since
there are at least $K-1$ other rows occupied by members in $T_\ell$ for
$\ell\neq k$.

As a result,
$\bar U_{i*}=\bar U_{j*}$ if $i,j\in T_{k}$ for some $k$,
and $\bar U_{i*}\neq\bar U_{j*}$ if $i\in T_{k}$, $j\in T_{\ell}$ with
$k\neq\ell$.
This gives a correspondence of clustering
between the rows in $\bar U_{T*}$ and those in $U_{T*}$ where $T=\bigcup_{k=1}^K T_k$.

\subsection{Proofs for degree corrected block models}\label{sec:pf-dcbm}

The argument fits very well in the
general argument developed in Section~\ref{sec:general}.
Then Lemma~\ref{lem:variational-bound} and Theorem~\ref{thm:long.thm}
still apply and
%
\begin{equation}
\label{eq:subspace-perturbation-adjacency}\qquad \pr \biggl[\|\hat U-UQ\|_F\le2\sqrt{2}C
\frac{\sqrt{Kn\alpha
_n}}{\gamma_n} \mbox{ for some } QQ^T=I_K \biggr]
\ge1-n^{-1} ,
\end{equation}
where $C$ is the constant in Theorem~\ref{thm:long.thm}.

For presentation simplicity, in the following argument we will work
with $Q=I_K$. The general
case can be handled in the same manner with more complicated notation
(simply substitute $U$
by $UQ$).

To prove Theorem~\ref{thm:main-dcbm}, we first give a bound on the
zero rows
in $\hat U$.
Recall that $I_+=\{i\dvtx\hat U_{i*}\neq0\}$. Define $I_0=I_+^c$.

\begin{lemma}[(Number of zero rows in $\hat U$)]\label{lem:discard-zero}
In a DCBM $(\Theta,B,\psi)$ satisfying the conditions of Theorem~\ref
{thm:main-dcbm}, let
$\hat U$ and $U$ be the leading eigenvectors of $A$ and~$P$, respectively.
Then
\[
|I_0|\le\sqrt{\sum_{k=1}^K
n_k^2 \nu_k}\|\hat U-U\|_F.
\]
\end{lemma}

\begin{pf} Use Cauchy--Schwarz:
\[
\|\hat U-U\|_F^2\ge\sum_{i=1}^n
\ind(\hat U_{i*}=0) \|U_{i*}\|^2\ge
\frac{ (\sum_{i=1}^n\ind(\hat U_{i*} = 0) )^{2}}{\sum_{i=1}^n
\|U_{i*}\|^{-2}}=\frac{|I_0|^{2}}{\sum_{k=1}^K
n_k^2 \nu_k}.
\]
%
\upqed\end{pf}

We also need the following simple fact about the distance between
normalized vectors.

\begin{fact*}
For two nonzero vectors $v_1$, $v_2$
of same dimension, we have
$\|\frac{v_1}{\|v_1\|}-\frac{v_2}{\|v_2\|}\|\le2\frac{\|v_1-v_2\|
}{\max
(\|v_1\|,\|v_2\|)}$.
\end{fact*}

\begin{pf}
Without loss of generality, assume $\|v_1\|\ge\|v_2\|$. Then
\begin{eqnarray*}
\biggl\llVert \frac{v_1}{\|v_1\|}-\frac{v_2}{\|v_2\|}\biggr\rrVert&=&\biggl\llVert
\frac{v_1}{\|v_1\|}-\frac{v_2}{\|v_1\|}+\frac{v_2}{\|v_1\|
}-\frac{v_2}{\|v_2\|}\biggr
\rrVert
\\
&\le& \frac{\|v_1-v_2\|}{\|v_1\|}+\frac{\|v_2\| |\|v_1\|-\|v_2\|
|}{\|v_1\|\|v_2\|} \le2\frac{\|v_1-v_2\|}{\|v_1\|}.
\end{eqnarray*}
\upqed\end{pf}

\begin{pf*}{Proof of Theorem~\ref{thm:main-dcbm}}
Recall that $U'$ is the row-normalized version of $U$.
Let $U''=U'_{I_+*}$ be the sub-matrix of $U'$ corresponding to
the nonzero rows in $\hat U$. Then
\begin{eqnarray*}
\bigl\|\hat U' - U''\bigr\|_{2,1}&\le&2 \sum
_{i=1}^n \frac{\|\hat U_{i*}-U_{i*}\|}{\|U_{i*}\|}
\\
&\le&2\sqrt{\sum_{i=1}^n \|\hat
U_{i*}-U_{i*}\|^2\sum
_{i=1}^n \|U_{i*}\|^{-2}}\le2
\sqrt{\|\hat U-U\|_F^2\sum
_{k=1}^K n_k^2
\nu_k}.
\end{eqnarray*}


Now we can bound the $(2,1)$ distance
between an approximate solution of $k$-median problem
\eqref{eq:k-median} and the targeted solution $U''$.
\begin{eqnarray*}
\bigl\|\hat\Theta^+ \hat X-U''\bigr\|_{2,1}&\le& \bigl\|\hat
\Theta^+ \hat X - \hat U'\bigr\|_{2,1}+ \bigl\|\hat U' -
U'' \bigr\|_{2,1}
\\
&\le& (2+\varepsilon)\bigl\|\hat U' - U''
\bigr\|_{2,1}.
\end{eqnarray*}

Let $S=\{i\in I_+\dvtx\|\hat\Theta_{i*} \hat X - U'_{i*}\|\ge\frac
{1}{\sqrt{2}}\}$.
The size of $S$
can be bounded using a similar argument as in the proof of Lemma~\ref
{lem:discard-zero}.
\begin{eqnarray*}
|S|\frac{1}{\sqrt{2}}&\le& \bigl\|\hat\Theta^+ \hat X-U''
\bigr\|_{2,1} \le (2+\varepsilon)\bigl\|\hat U' - U''
\bigr\|_{2,1}
\\
&\le& 2(2+\varepsilon)\sqrt{\sum_{k=1}^K
n_k^2 \nu_k}\|\hat U-U\|_F
,
\end{eqnarray*}
which implies
%
\begin{equation}
\label{eq:dcbm-bound-S} |S|\le 2\sqrt{2}(2+\varepsilon)\sqrt{\sum
_{k=1}^K n_k^2
\nu_k}\|\hat U-U\| _F.%
\end{equation}

On the event in \eqref{eq:subspace-perturbation-adjacency} (recall that
we assume $Q=I$),
\eqref{eq:dcbm-bound-S} and Lemma~\ref{lem:discard-zero} implies
%
\begin{equation}
\label{eq:dcbm-bound-S-zero} |S|+|I_0| \le (2.5+\varepsilon)8C\frac{\sqrt{K n \alpha_n}}{\gamma_n}
\sqrt {\sum_{k=1}^K
n_k^2\nu_k}.
\end{equation}
Combining this with condition
\eqref{eq:dcbm-main-condition} implies
$|S|+|I_0|< n_k$ for all $k$ and hence
$G_k\cap(I_+ \setminus S)\neq\varnothing$.
Therefore, for any two rows in $G:= I_+\setminus S$,
if they are in different clusters of $\Theta$ then they must be in
different clusters of $\hat\Theta$
(otherwise, $\|U_{i*}'-U_{j*}'\|\le\|U_{i*}'-\hat\Theta_{i*} \hat X\|
+\|\hat\Theta_{j*}\hat X-U_{j*}'\|
< \sqrt{2}$).

As a consequence, the mis-clustered nodes are no more than $I_0\cup S$,
and the number
is bounded by the right-hand side of \eqref{eq:dcbm-bound-S-zero}. The
claimed result follows
by choosing $c=8C$.
\end{pf*}
\end{appendix}


\section*{Acknowledgment}
The authors thank an anonymous reviewer for helpful suggestions that
led in
particular to a
significant simplification of the proof of
Lemma~\ref{lem:variational-bound}.


\begin{supplement}[id=suppA]
\stitle{Supplement to ``Consistency of spectral clustering in sparse
stochastic block models''}
\slink[doi]{10.1214/14-AOS1274SUPP} 
\sdatatype{.pdf}
\sfilename{aos1274\_supp.pdf}
\sdescription{The supplementary file contains a proof of Theorem~\ref
{thm:long.thm}.}
\end{supplement}

%
%
%

\printaddresses

\begin{thebibliography}{48}

\bibitem[\protect\citeauthoryear{Aloise et~al.}{2009}]{AloiseDHP09}
%
\begin{barticle}[author]
\bauthor{\bsnm{Aloise},~\bfnm{Daniel}\binits{D.}},
\bauthor{\bsnm{Deshpande},~\bfnm{Amit}\binits{A.}},
\bauthor{\bsnm{Hansen},~\bfnm{Pierre}\binits{P.}} \AND
\bauthor{\bsnm{Popat},~\bfnm{Preyas}\binits{P.}}
(\byear{2009}).
\btitle{NP-hardness of Euclidean sum-of-squares clustering}.
\bjournal{Machine Learning}
\bvolume{75}
\bpages{245--248}.
\end{barticle}
%
\bptok{imsref}%
\endbibitem

\bibitem[\protect\citeauthoryear{Alon and Spencer}{2004}]{AlonS04}
%
\begin{bbook}[author]
\bauthor{\bsnm{Alon},~\bfnm{Noga}\binits{N.}} \AND
\bauthor{\bsnm{Spencer},~\bfnm{Joel~H.}\binits{J.~H.}}
(\byear{2004}).
\btitle{The Probabilistic Method},
\bedition{2nd} ed.
\bpublisher{Wiley},
\blocation{Hoboken}.
\end{bbook}
\bptok{imsref}%
\endbibitem

\bibitem[\protect\citeauthoryear{Amini et~al.}{2012}]{AminiCBL12}
%
\begin{bmisc}[author]
\bauthor{\bsnm{Amini},~\bfnm{Arash~A.}\binits{A.~A.}},
\bauthor{\bsnm{Chen},~\bfnm{Aiyou}\binits{A.}},
\bauthor{\bsnm{Bickel},~\bfnm{Peter~J.}\binits{P.~J.}} \AND
\bauthor{\bsnm{Levina},~\bfnm{Elizaveta}\binits{E.}}
(\byear{2012}).
\bhowpublished{Pseudo-likelihood methods for community detection
in large sparse networks.
Preprint. Available at \arxivurl{arXiv:1207.2340}}.
\end{bmisc}
%
\bptok{imsref}%
\endbibitem

\bibitem[\protect\citeauthoryear{Anandkumar et~al.}{2013}]{Anandkumar13}
%
\begin{bmisc}[author]
\bauthor{\bsnm{Anandkumar},~\bfnm{Anima}\binits{A.}},
\bauthor{\bsnm{Ge},~\bfnm{Rong}\binits{R.}},
\bauthor{\bsnm{Hsu},~\bfnm{Daniel}\binits{D.}} \AND
\bauthor{\bsnm{Kakade},~\bfnm{Sham~M.}\binits{S.~M.}}
(\byear{2013}).
\bhowpublished{A tensor spectral approach to learning mixed membership
community models.
Preprint. Available at \arxivurl{arXiv:1302.2684}}.
\end{bmisc}
%
\bptok{imsref}%
\endbibitem

\bibitem[\protect\citeauthoryear{Awasthi and Sheffet}{2012}]{AwasthiS12}
%
\begin{bincollection}[mr]
\bauthor{\bsnm{Awasthi},~\bfnm{Pranjal}\binits{P.}} \AND
\bauthor{\bsnm{Sheffet},~\bfnm{Or}\binits{O.}}
(\byear{2012}).
\btitle{Improved spectral-norm bounds for clustering}.
In \bbooktitle{Approximation, Randomization, and Combinatorial Optimization}.
\bseries{Lecture Notes in Computer Science}
\bvolume{7408}
\bpages{37--49}.
\bpublisher{Springer},
\blocation{Heidelberg}.
\bid{doi={10.1007/978-3-642-32512-0_4}, mr={3003539}}
\end{bincollection}
%
\bptok{imsref}%
\endbibitem

\bibitem[\protect\citeauthoryear{Balakrishnan
et~al.}{2011}]{siva.spectral:11}
%
\begin{binproceedings}[author]
\bauthor{\bsnm{Balakrishnan},~\bfnm{Sivaraman}\binits{S.}},
\bauthor{\bsnm{Xu},~\bfnm{Min}\binits{M.}},
\bauthor{\bsnm{Krishnamurthy},~\bfnm{Akshay}\binits{A.}} \AND
\bauthor{\bsnm{Singh},~\bfnm{Aarti}\binits{A.}}
(\byear{2011}).
\btitle{Noise thresholds for spectral clustering}.
In \bbooktitle{Advances in Neural Information Processing Systems 24}
(\beditor{\binits{J.}\bfnm{J.} \bsnm{Shawe-Taylor}},
\beditor{\binits{R.~S.}\bfnm{R.~S.} \bsnm{Zemel}},
\beditor{\binits{P.~L.}\bfnm{P.~L.} \bsnm{Bartlett}},
\beditor{\binits{F.}\bfnm{F.} \bsnm{Pereira}}
\AND
\beditor{\binits{K.~Q.}\bfnm{K.~Q.} \bsnm{Weinberger}}, eds.)
\bpages{954--962}.
\bpublisher{Curran Associates},
\blocation{Red Hook, NY}.
\end{binproceedings}
%
\bptok{imsref}%
\endbibitem

\bibitem[\protect\citeauthoryear{Bhatia}{1997}]{Bhatia97}
%
\begin{bbook}[mr]
\bauthor{\bsnm{Bhatia},~\bfnm{Rajendra}\binits{R.}}
(\byear{1997}).
\btitle{Matrix Analysis}.
\bseries{Graduate Texts in Mathematics}
\bvolume{169}.
\bpublisher{Springer},
\blocation{New York}.
\bid{doi={10.1007/978-1-4612-0653-8}, mr={1477662}}
\end{bbook}
%
\bptok{imsref}%
\endbibitem

\bibitem[\protect\citeauthoryear{Bickel and Chen}{2009}]{BickelC09}
%
\begin{barticle}[author]
\bauthor{\bsnm{Bickel},~\bfnm{Peter~J.}\binits{P.~J.}} \AND
\bauthor{\bsnm{Chen},~\bfnm{Aiyou}\binits{A.}}
(\byear{2009}).
\btitle{A nonparametric view of network models and Newman--Girvan and
other modularities}.
\bjournal{Proc. Natl. Acad. Sci. USA}
\bvolume{106}
\bpages{21068--21073}.
\end{barticle}
%
\bptok{imsref}%
\endbibitem

\bibitem[\protect\citeauthoryear{Celisse, Daudin and
Pierre}{2012}]{Celisse12}
%
\begin{barticle}[mr]
\bauthor{\bsnm{Celisse},~\bfnm{Alain}\binits{A.}},
\bauthor{\bsnm{Daudin},~\bfnm{Jean-Jacques}\binits{J.-J.}} \AND
\bauthor{\bsnm{Pierre},~\bfnm{Laurent}\binits{L.}}
(\byear{2012}).
\btitle{Consistency of maximum-likelihood and variational estimators
in the stochastic block model}.
\bjournal{Electron. J. Stat.}
\bvolume{6}
\bpages{1847--1899}.
\bid{doi={10.1214/12-EJS729}, issn={1935-7524}, mr={2988467}}
\end{barticle}
%
\bptok{imsref}%
\endbibitem

\bibitem[\protect\citeauthoryear{Channarond, Daudin and
Robin}{2012}]{stephane}
%
\begin{barticle}[mr]
\bauthor{\bsnm{Channarond},~\bfnm{Antoine}\binits{A.}},
\bauthor{\bsnm{Daudin},~\bfnm{Jean-Jacques}\binits{J.-J.}} \AND
\bauthor{\bsnm{Robin},~\bfnm{St{\'e}phane}\binits{S.}}
(\byear{2012}).
\btitle{Classification and estimation in the stochastic blockmodel
based on the empirical degrees}.
\bjournal{Electron. J. Stat.}
\bvolume{6}
\bpages{2574--2601}.
\bid{doi={10.1214/12-EJS753}, issn={1935-7524}, mr={3020277}}
\end{barticle}
%
\bptok{imsref}%
\endbibitem

\bibitem[\protect\citeauthoryear{Charikar et~al.}{1999}]{CharikarGTS99}
%
\begin{binproceedings}[author]
\bauthor{\bsnm{Charikar},~\bfnm{Moses}\binits{M.}},
\bauthor{\bsnm{Guha},~\bfnm{Sudipto}\binits{S.}},
\bauthor{\bsnm{Tardos},~\bfnm{{\'E}va}\binits{{\'E}.}} \AND
\bauthor{\bsnm{Shmoys},~\bfnm{David~B.}\binits{D.~B.}}
(\byear{1999}).
\btitle{A constant-factor approximation algorithm for the $k$-median problem}.
In \bbooktitle{Proceedings of the Thirty-First Annual ACM Symposium on
Theory of Computing}
\bpages{1--10}.
\bpublisher{ACM},
\blocation{New York, NY}.
\end{binproceedings}
%
\bptok{imsref}%
\endbibitem

\bibitem[\protect\citeauthoryear{Chaudhuri, Chung and
Tsiatas}{2012}]{ChaudhuriCT12}
%
\begin{barticle}[author]
\bauthor{\bsnm{Chaudhuri},~\bfnm{Kamalika}\binits{K.}},
\bauthor{\bsnm{Chung},~\bfnm{Fan}\binits{F.}} \AND
\bauthor{\bsnm{Tsiatas},~\bfnm{Alexander}\binits{A.}}
(\byear{2012}).
\btitle{Spectral clustering of graphs with general degrees in the
extended planted partition model}.
\bjournal{JMLR: Workshop and Conference Proceedings}
\bvolume{2012}
\bpages{35.1--35.23}.
\end{barticle}
%
\bptok{imsref}%
\endbibitem

\bibitem[\protect\citeauthoryear{Chen, Sanghavi and Xu}{2012}]{ChenSX12}
%
\begin{bincollection}[author]
\bauthor{\bsnm{Chen},~\bfnm{Yudong}\binits{Y.}},
\bauthor{\bsnm{Sanghavi},~\bfnm{Sujay}\binits{S.}} \AND
\bauthor{\bsnm{Xu},~\bfnm{Huan}\binits{H.}}
(\byear{2012}).
\btitle{Clustering sparse graphs}.
In \bbooktitle{Advances in Neural Information Processing Systems 25}
(\beditor{\bfnm{F.}\binits{F.}~\bsnm{Pereira}},
\beditor{\bfnm{C.~J.~C.}\binits{C.~J.~C.}~\bsnm{Burges}},
\beditor{\bfnm{L.}\binits{L.}~\bsnm{Bottou}} \AND
\beditor{\bfnm{K.~Q.}\binits{K.~Q.}~\bsnm{Weinberger}}, eds.)
\bpages{2204--2212}.
\bpublisher{Curran Associates},
\blocation{Red Hook, NY}.
\end{bincollection}
%
\bptok{imsref}%
\endbibitem

\bibitem[\protect\citeauthoryear{Choi, Wolfe and Airoldi}{2012}]{ChoiWA12}
%
\begin{barticle}[mr]
\bauthor{\bsnm{Choi},~\bfnm{D.~S.}\binits{D.~S.}},
\bauthor{\bsnm{Wolfe},~\bfnm{P.~J.}\binits{P.~J.}} \AND
\bauthor{\bsnm{Airoldi},~\bfnm{E.~M.}\binits{E.~M.}}
(\byear{2012}).
\btitle{Stochastic blockmodels with a growing number of classes}.
\bjournal{Biometrika}
\bvolume{99}
\bpages{273--284}.
\bid{doi={10.1093/biomet/asr053}, issn={0006-3444}, mr={2931253}}
\end{barticle}
%
\bptok{imsref}%
\endbibitem

\bibitem[\protect\citeauthoryear{Chung and Radcliffe}{2011}]{ChungR11}
%
\begin{barticle}[mr]
\bauthor{\bsnm{Chung},~\bfnm{Fan}\binits{F.}} \AND
\bauthor{\bsnm{Radcliffe},~\bfnm{Mary}\binits{M.}}
(\byear{2011}).
\btitle{On the spectra of general random graphs}.
\bjournal{Electron. J. Combin.}
\bvolume{18}
\bpages{Paper 215, 14}.
\bid{issn={1077-8926}, mr={2853072}}
\end{barticle}
%
\bptok{imsref}%
\endbibitem

\bibitem[\protect\citeauthoryear{Coja-Oghlan}{2010}]{Coja-Oghlan10}
%
\begin{barticle}[mr]
\bauthor{\bsnm{Coja-Oghlan},~\bfnm{Amin}\binits{A.}}
(\byear{2010}).
\btitle{Graph partitioning via adaptive spectral techniques}.
\bjournal{Combin. Probab. Comput.}
\bvolume{19}
\bpages{227--284}.
\bid{doi={10.1017/S0963548309990514}, issn={0963-5483}, mr={2593622}}
\end{barticle}
%
\bptok{imsref}%
\endbibitem

\bibitem[\protect\citeauthoryear{Decelle et~al.}{2011}]{Decelle11}
%
\begin{barticle}[author]
\bauthor{\bsnm{Decelle},~\bfnm{Aurelien}\binits{A.}},
\bauthor{\bsnm{Krzakala},~\bfnm{Florent}\binits{F.}},
\bauthor{\bsnm{Moore},~\bfnm{Cristopher}\binits{C.}} \AND
\bauthor{\bsnm{Zdeborov{\'a}},~\bfnm{Lenka}\binits{L.}}
(\byear{2011}).
\btitle{Asymptotic analysis of the stochastic block model for modular
networks and its algorithmic applications}.
\bjournal{Phys. Rev. E~(3)}
\bvolume{84}
\bpages{066106}.
\end{barticle}
%
\bptok{imsref}%
\endbibitem

\bibitem[\protect\citeauthoryear{Deshpande and
Montanari}{2013}]{DeshpandeM13}
%
\begin{bmisc}[author]
\bauthor{\bsnm{Deshpande},~\bfnm{Yash}\binits{Y.}} \AND
\bauthor{\bsnm{Montanari},~\bfnm{Andrea}\binits{A.}}
(\byear{2013}).
\bhowpublished{Finding hidden cliques of size $\sqrt{N/e}$ in nearly
linear time.
Preprint.
Available at \arxivurl{arXiv:1304.7047}}.
\end{bmisc}
%
\bptok{imsref}%
\endbibitem

\bibitem[\protect\citeauthoryear{Feige and Ofek}{2005}]{FeigeO05}
%
\begin{barticle}[mr]
\bauthor{\bsnm{Feige},~\bfnm{Uriel}\binits{U.}} \AND
\bauthor{\bsnm{Ofek},~\bfnm{Eran}\binits{E.}}
(\byear{2005}).
\btitle{Spectral techniques applied to sparse random graphs}.
\bjournal{Random Structures Algorithms}
\bvolume{27}
\bpages{251--275}.
\bid{doi={10.1002/rsa.20089}, issn={1042-9832}, mr={2155709}}
\end{barticle}
%
\bptok{imsref}%
\endbibitem

\bibitem[\protect\citeauthoryear{Fishkind et~al.}{2013}]{Fishkind13}
%
\begin{barticle}[mr]
\bauthor{\bsnm{Fishkind},~\bfnm{Donniell~E.}\binits{D.~E.}},
\bauthor{\bsnm{Sussman},~\bfnm{Daniel~L.}\binits{D.~L.}},
\bauthor{\bsnm{Tang},~\bfnm{Minh}\binits{M.}},
\bauthor{\bsnm{Vogelstein},~\bfnm{Joshua~T.}\binits{J.~T.}} \AND
\bauthor{\bsnm{Priebe},~\bfnm{Carey~E.}\binits{C.~E.}}
(\byear{2013}).
\btitle{Consistent adjacency-spectral partitioning for the stochastic
block model when the model parameters are unknown}.
\bjournal{SIAM J. Matrix Anal. Appl.}
\bvolume{34}
\bpages{23--39}.
\bid{doi={10.1137/120875600}, issn={0895-4798}, mr={3032990}}
\end{barticle}
%
\bptok{imsref}%
\endbibitem

\bibitem[\protect\citeauthoryear{Goldenberg et~al.}{2010}]{Goldenberg10}
%
\begin{barticle}[author]
\bauthor{\bsnm{Goldenberg},~\bfnm{Anna}\binits{A.}},
\bauthor{\bsnm{Zheng},~\bfnm{Alice~X.}\binits{A.~X.}},
\bauthor{\bsnm{Fienberg},~\bfnm{Stephen~E.}\binits{S.~E.}} \AND
\bauthor{\bsnm{Airoldi},~\bfnm{Edoardo~M.}\binits{E.~M.}}
(\byear{2010}).
\btitle{A survey of statistical network models}.
\bjournal{Foundations and Trends{\textregistered} in Machine Learning}
\bvolume{2}
\bpages{129--233}.
\end{barticle}
%
\bptok{imsref}%
\endbibitem

\bibitem[\protect\citeauthoryear{Holland, Laskey and
Leinhardt}{1983}]{Holland83}
%
\begin{barticle}[mr]
\bauthor{\bsnm{Holland},~\bfnm{Paul~W.}\binits{P.~W.}},
\bauthor{\bsnm{Laskey},~\bfnm{Kathryn~Blackmond}\binits{K.~B.}}
\AND
\bauthor{\bsnm{Leinhardt},~\bfnm{Samuel}\binits{S.}}
(\byear{1983}).
\btitle{Stochastic blockmodels: First steps}.
\bjournal{Social Networks}
\bvolume{5}
\bpages{109--137}.
\bid{doi={10.1016/0378-8733(83)90021-7}, issn={0378-8733}, mr={0718088}}
\end{barticle}
%
\bptok{imsref}%
\endbibitem

\bibitem[\protect\citeauthoryear{Jin}{2012}]{Jin12}
%
\begin{bmisc}[author]
\bauthor{\bsnm{Jin},~\bfnm{Jiashun}\binits{J.}}
(\byear{2012}).
\bhowpublished{Fast community detection by SCORE.
Available at \arxivurl{arXiv:1211.5803}}.
\end{bmisc}
%
\bptok{imsref}%
\endbibitem

\bibitem[\protect\citeauthoryear{Karrer and Newman}{2011}]{KarrerN11}
%
\begin{barticle}[mr]
\bauthor{\bsnm{Karrer},~\bfnm{Brian}\binits{B.}} \AND
\bauthor{\bsnm{Newman},~\bfnm{M.~E.~J.}\binits{M.~E.~J.}}
(\byear{2011}).
\btitle{Stochastic blockmodels and community structure in networks}.
\bjournal{Phys. Rev. E (3)}
\bvolume{83}
\bpages{016107, 10}.
\bid{doi={10.1103/PhysRevE.83.016107}, issn={1539-3755}, mr={2788206}}
\end{barticle}
%
\bptok{imsref}%
\endbibitem

\bibitem[\protect\citeauthoryear{Kolaczyk}{2009}]{Kolaczyk09}
%
\begin{bbook}[mr]
\bauthor{\bsnm{Kolaczyk},~\bfnm{Eric~D.}\binits{E.~D.}}
(\byear{2009}).
\btitle{Statistical Analysis of Network Data: Methods and Models}.
\bpublisher{Springer},
\blocation{New York}.
\bid{doi={10.1007/978-0-387-88146-1}, mr={2724362}}
\end{bbook}
%
\bptok{imsref}%
\endbibitem

\bibitem[\protect\citeauthoryear{Krzakala et~al.}{2013}]{Krzakala13}
%
\begin{barticle}[mr]
\bauthor{\bsnm{Krzakala},~\bfnm{Florent}\binits{F.}},
\bauthor{\bsnm{Moore},~\bfnm{Cristopher}\binits{C.}},
\bauthor{\bsnm{Mossel},~\bfnm{Elchanan}\binits{E.}},
\bauthor{\bsnm{Neeman},~\bfnm{Joe}\binits{J.}},
\bauthor{\bsnm{Sly},~\bfnm{Allan}\binits{A.}},
\bauthor{\bsnm{Zdeborov{\'a}},~\bfnm{Lenka}\binits{L.}} \AND
\bauthor{\bsnm{Zhang},~\bfnm{Pan}\binits{P.}}
(\byear{2013}).
\btitle{Spectral redemption in clustering sparse networks}.
\bjournal{Proc. Natl. Acad. Sci. USA}
\bvolume{110}
\bpages{20935--20940}.
\bid{doi={10.1073/pnas.1312486110}, issn={1091-6490}, mr={3174850}}
\end{barticle}
%
\bptok{imsref}%
\endbibitem

\bibitem[\protect\citeauthoryear{Kumar and Kannan}{2010}]{KumarK10}
%
\begin{bincollection}[mr]
\bauthor{\bsnm{Kumar},~\bfnm{Amit}\binits{A.}} \AND
\bauthor{\bsnm{Kannan},~\bfnm{Ravindran}\binits{R.}}
(\byear{2010}).
\btitle{Clustering with spectral norm and the {$k$}-means algorithm}.
In \bbooktitle{Proceedings of the 2010 IEEE 51st {A}nnual {S}ymposium on {F}oundations of
{C}omputer {S}cience FOCS}
\bpages{299--308}.
\bpublisher{IEEE},
\blocation{Los Alamitos, CA}.
\bid{mr={3025203}}
\end{bincollection}
%
\bptok{imsref}%
\endbibitem

\bibitem[\protect\citeauthoryear{Kumar, Sabharwal and Sen}{2004}]{KumarSS04}
%
\begin{binproceedings}[author]
\bauthor{\bsnm{Kumar},~\bfnm{Amit}\binits{A.}},
\bauthor{\bsnm{Sabharwal},~\bfnm{Yogish}\binits{Y.}} \AND
\bauthor{\bsnm{Sen},~\bfnm{Sandeep}\binits{S.}}
(\byear{2004}).
\btitle{A simple linear time ($1+\varepsilon$)-approximation
algorithm for $k$-means clustering in any dimensions}.
In \bbooktitle{Proceedings of the
45th Annual IEEE Symposium on
Foundations of Computer Science}
\bpages{454--462}.
\bpublisher{IEEE Computer Society},
\blocation{Washington, DC}.
\end{binproceedings}
%
\bptok{imsref}%
\endbibitem

\bibitem[\protect\citeauthoryear{Lei and Rinaldo}{2014}]{SUPP}
%
\begin{bmisc}[author]
\bauthor{\bsnm{Lei},~\bfnm{Jing}\binits{J.}} \AND
\bauthor{\bsnm{Rinaldo},~\bfnm{Alessandro}\binits{A.}}
(\byear{2014}).
\bhowpublished{Supplement to ``Consistency of spectral clustering in
stochastic block models.''
DOI:\doiurl{10.1214/14-AOS1274SUPP}}.
\bptok{imsref}%
\end{bmisc}
%
\endbibitem
%
%
%
\bptok{imsref}%
\endbibitem

\bibitem[\protect\citeauthoryear{Li and Svensson}{2013}]{LiS13}
%
\begin{binproceedings}[author]
\bauthor{\bsnm{Li},~\bfnm{Shi}\binits{S.}} \AND
\bauthor{\bsnm{Svensson},~\bfnm{Ola}\binits{O.}}
(\byear{2013}).
\btitle{Approximating k-median via pseudo-approximation}.
In \bbooktitle{Proceedings of the 45th Annual ACM Symposium on
Symposium on Theory of Computing}
\bpages{901--910}.
\bpublisher{ACM},
\blocation{New York}.
\end{binproceedings}
%
\bptok{imsref}%
\endbibitem

\bibitem[\protect\citeauthoryear{Lu and Peng}{2012}]{LuP12}
%
\begin{bmisc}[author]
\bauthor{\bsnm{Lu},~\bfnm{Linyuan}\binits{L.}} \AND
\bauthor{\bsnm{Peng},~\bfnm{Xing}\binits{X.}}
(\byear{2012}).
\bhowpublished{Spectra of edge-independent random graphs.
Preprint. Available at \arxivurl{arXiv:1204.6207}}.
\end{bmisc}
%
\bptok{imsref}%
\endbibitem

\bibitem[\protect\citeauthoryear{Lyzinski et~al.}{2013}]{Lyzinski13}
%
\begin{bmisc}[author]
\bauthor{\bsnm{Lyzinski},~\bfnm{Vince}\binits{V.}},
\bauthor{\bsnm{Sussman},~\bfnm{Daniel}\binits{D.}},
\bauthor{\bsnm{Tang},~\bfnm{Minh}\binits{M.}},
\bauthor{\bsnm{Athreya},~\bfnm{Avanti}\binits{A.}} \AND
\bauthor{\bsnm{Priebe},~\bfnm{Carey}\binits{C.}}
(\byear{2013}).
\bhowpublished{Perfect clustering for stochastic blockmodel graphs
via adjacency spectral embedding.
Preprint. Available at \arxivurl{arXiv:1310.0532}.}
\end{bmisc}
%
\bptok{imsref}%
\endbibitem

\bibitem[\protect\citeauthoryear{Massoulie}{2013}]{Massoulie13}
%
\begin{bmisc}[author]
\bauthor{\bsnm{Massoulie},~\bfnm{Laurent}\binits{L.}}
(\byear{2013}).
\bhowpublished{Community detection thresholds and the weak Ramanujan property.
Preprint. Available at \arxivurl{arXiv:1311.3085}}.
\end{bmisc}
%
\bptok{imsref}%
\endbibitem

\bibitem[\protect\citeauthoryear{McSherry}{2001}]{McSherry01}
%
\begin{bincollection}[mr]
\bauthor{\bsnm{McSherry},~\bfnm{Frank}\binits{F.}}
(\byear{2001}).
\btitle{Spectral partitioning of random graphs}.
In \bbooktitle{42nd IEEE {S}ymposium on {F}oundations of {C}omputer
{S}cience ({L}as {V}egas, NV, 2001)}
\bpages{529--537}.
\bpublisher{IEEE},
\blocation{Los Alamitos, CA}.
\bid{mr={1948742}}
\end{bincollection}
\bptok{imsref}%
\endbibitem

\bibitem[\protect\citeauthoryear{Mossel, Neeman and Sly}{2012}]{MosselNS12}
%
\begin{bmisc}[author]
\bauthor{\bsnm{Mossel},~\bfnm{Elchanan}\binits{E.}},
\bauthor{\bsnm{Neeman},~\bfnm{Joe}\binits{J.}} \AND
\bauthor{\bsnm{Sly},~\bfnm{Allan}\binits{A.}}
(\byear{2012}).
\bhowpublished{Stochastic block models and reconstruction.
Preprint. Available at \arxivurl{arXiv:1202.1499}}.
\end{bmisc}
%
\bptok{imsref}%
\endbibitem

\bibitem[\protect\citeauthoryear{Mossel, Neeman and Sly}{2013}]{MosselNS13}
%
\begin{bmisc}[author]
\bauthor{\bsnm{Mossel},~\bfnm{Elchanan}\binits{E.}},
\bauthor{\bsnm{Neeman},~\bfnm{Joe}\binits{J.}} \AND
\bauthor{\bsnm{Sly},~\bfnm{Allan}\binits{A.}}
(\byear{2013}).
\bhowpublished{A proof of the block model threshold conjecture.
Preprint. Available at \arxivurl{arXiv:1311.4115}}.
\end{bmisc}
%
\bptok{imsref}%
\endbibitem

\bibitem[\protect\citeauthoryear{Newman}{2010}]{Newman09}
%
\begin{bbook}[mr]
\bauthor{\bsnm{Newman},~\bfnm{M.~E.~J.}\binits{M.~E.~J.}}
(\byear{2010}).
\btitle{Networks: An Introduction}.
\bpublisher{Oxford Univ. Press},
\blocation{Oxford}.
\bid{doi={10.1093/acprof:oso/9780199206650.001.0001}, mr={2676073}}
\bptnote{check year}%
\end{bbook}
%
\bptok{imsref}%
\endbibitem

\bibitem[\protect\citeauthoryear{Newman and Girvan}{2004}]{NewmanG04}
%
\begin{barticle}[author]
\bauthor{\bsnm{Newman},~\bfnm{Mark~E.~J.}\binits{M.~E.~J.}} \AND
\bauthor{\bsnm{Girvan},~\bfnm{Michelle}\binits{M.}}
(\byear{2004}).
\btitle{Finding and evaluating community structure in networks}.
\bjournal{Phys. Rev. E (3)}
\bvolume{69}
\bpages{026113}.
\end{barticle}
%
\bptok{imsref}%
\endbibitem

\bibitem[\protect\citeauthoryear{Ng et~al.}{2002}]{NgJW02}
%
\begin{barticle}[author]
\bauthor{\bsnm{Ng},~\bfnm{Andrew~Y.}\binits{A.~Y.}},
\bauthor{\bsnm{Jordan},~\bfnm{Michael~I.}\binits{M.~I.}},
\bauthor{\bsnm{Weiss},~\bfnm{Yair}\binits{Y.}} \betal{et~al.}
(\byear{2002}).
\btitle{On spectral clustering: Analysis and an algorithm}.
\bjournal{Adv. Neural Inf. Process. Syst.}
\bvolume{2}
\bpages{849--856}.
\end{barticle}
%
\bptok{imsref}%
\endbibitem

\bibitem[\protect\citeauthoryear{Qin and Rohe}{2013}]{QinR13}
%
\begin{bmisc}[author]
\bauthor{\bsnm{Qin},~\bfnm{Tai}\binits{T.}} \AND
\bauthor{\bsnm{Rohe},~\bfnm{Karl}\binits{K.}}
(\byear{2013}).
\bhowpublished{Regularized spectral clustering under the
degree-corrected stochastic blockmodel.
Preprint. Available at
\arxivurl{arXiv:1309.4111}}.
\end{bmisc}
%
\bptok{imsref}%
\endbibitem

\bibitem[\protect\citeauthoryear{Rohe, Chatterjee and Yu}{2011}]{RoheCY11}
%
\begin{barticle}[mr]
\bauthor{\bsnm{Rohe},~\bfnm{Karl}\binits{K.}},
\bauthor{\bsnm{Chatterjee},~\bfnm{Sourav}\binits{S.}} \AND
\bauthor{\bsnm{Yu},~\bfnm{Bin}\binits{B.}}
(\byear{2011}).
\btitle{Spectral clustering and the high-dimensional stochastic blockmodel}.
\bjournal{Ann. Statist.}
\bvolume{39}
\bpages{1878--1915}.
\bid{doi={10.1214/11-AOS887}, issn={0090-5364}, mr={2893856}}
\end{barticle}
%
\bptok{imsref}%
\endbibitem

\bibitem[\protect\citeauthoryear{Sarkar and Bickel}{2013}]{sarkar.bickel:13}
%
\begin{bmisc}[author]
\bauthor{\bsnm{Sarkar},~\bfnm{Purnamrita}\binits{P.}} \AND
\bauthor{\bsnm{Bickel},~\bfnm{Peter}\binits{P.}}
(\byear{2013}).
\bhowpublished{Role of normalization in spectral clustering
for stochastic blockmodels.
Preprint. Available at \arxivurl{arXiv:1310.1495}}.
\end{bmisc}
%
\bptok{imsref}%
\endbibitem

\bibitem[\protect\citeauthoryear{Sussman et~al.}{2012}]{Sussman:12}
%
\begin{barticle}[mr]
\bauthor{\bsnm{Sussman},~\bfnm{Daniel~L.}\binits{D.~L.}},
\bauthor{\bsnm{Tang},~\bfnm{Minh}\binits{M.}},
\bauthor{\bsnm{Fishkind},~\bfnm{Donniell~E.}\binits{D.~E.}} \AND
\bauthor{\bsnm{Priebe},~\bfnm{Carey~E.}\binits{C.~E.}}
(\byear{2012}).
\btitle{A consistent adjacency spectral embedding for stochastic
blockmodel graphs}.
\bjournal{J. Amer. Statist. Assoc.}
\bvolume{107}
\bpages{1119--1128}.
\bid{doi={10.1080/01621459.2012.699795}, issn={0162-1459}, mr={3010899}}
\end{barticle}
%
\bptok{imsref}%
\endbibitem


\bibitem[\protect\citeauthoryear{Tropp}{2012}]{Tropp12}
%
\begin{barticle}[mr]
\bauthor{\bsnm{Tropp},~\bfnm{Joel~A.}\binits{J.~A.}}
(\byear{2012}).
\btitle{User-friendly tail bounds for sums of random matrices}.
\bjournal{Found. Comput. Math.}
\bvolume{12}
\bpages{389--434}.
\bid{doi={10.1007/s10208-011-9099-z}, issn={1615-3375}, mr={2946459}}
\end{barticle}
%
\bptok{imsref}%
\endbibitem

\bibitem[\protect\citeauthoryear{von Luxburg}{2007}]{VonLuxburg07}
%
\begin{barticle}[mr]
\bauthor{\bparticle{von} \bsnm{Luxburg},~\bfnm{Ulrike}\binits{U.}}
(\byear{2007}).
\btitle{A tutorial on spectral clustering}.
\bjournal{Stat. Comput.}
\bvolume{17}
\bpages{395--416}.
\bid{doi={10.1007/s11222-007-9033-z}, issn={0960-3174}, mr={2409803}}
\end{barticle}
%
\bptok{imsref}%
\endbibitem

\bibitem[\protect\citeauthoryear{Vu and Lei}{2013}]{VuL13}
%
\begin{barticle}[mr]
\bauthor{\bsnm{Vu},~\bfnm{Vincent~Q.}\binits{V.~Q.}} \AND
\bauthor{\bsnm{Lei},~\bfnm{Jing}\binits{J.}}
(\byear{2013}).
\btitle{Minimax sparse principal subspace estimation in high dimensions}.
\bjournal{Ann. Statist.}
\bvolume{41}
\bpages{2905--2947}.
\bid{doi={10.1214/13-AOS1151}, issn={0090-5364}, mr={3161452}}
\end{barticle}
%
\bptok{imsref}%
\endbibitem

\bibitem[\protect\citeauthoryear{Zhao, Levina and Zhu}{2012}]{ZhaoLZ12}
%
\begin{barticle}[mr]
\bauthor{\bsnm{Zhao},~\bfnm{Yunpeng}\binits{Y.}},
\bauthor{\bsnm{Levina},~\bfnm{Elizaveta}\binits{E.}} \AND
\bauthor{\bsnm{Zhu},~\bfnm{Ji}\binits{J.}}
(\byear{2012}).
\btitle{Consistency of community detection in networks under
degree-corrected stochastic block models}.
\bjournal{Ann. Statist.}
\bvolume{40}
\bpages{2266--2292}.
\bid{doi={10.1214/12-AOS1036}, issn={0090-5364}, mr={3059083}}
\end{barticle}
%
\bptok{imsref}%
\endbibitem
\end{thebibliography}
\end{document}